\documentclass{article}

\usepackage{amsmath,amssymb,latexsym,amsthm,epsfig,stmaryrd,tocbibind}
\usepackage[OT2,OT1]{fontenc}

\newcommand\cyr{%
\renewcommand\rmdefault{wncyr}%
\renewcommand\sfdefault{wncyss}%
\renewcommand\encodingdefault{OT2}%
\normalfont \selectfont} \DeclareTextFontCommand{\textcyr}{\cyr}
\def\part#1{\frac{\partial\phantom{q}}{\partial#1}}

\newenvironment{rmk}{\begin{trivlist}\item[]{\bf Remark:} }
{\end{trivlist}}
\newenvironment{ex}{\begin{trivlist}\item[]{\bf Example:} }
{\end{trivlist}}
\newenvironment{prf}{\begin{trivlist}\item[]{\bf Proof:} }
{\hfill $\blacksquare$ \end{trivlist}}

\newtheorem{thm}{Theorem}[section]
\newtheorem{definition}{Definition}[section]
\newtheorem{prp}[thm]{Proposition}
\newtheorem{lem}[thm]{Lemma}
\newtheorem{cor}[thm]{Corollary}

\def\End{\mathop{\rm End}\nolimits}

\def\ker{\mathop{\rm ker}\nolimits}

\def\cliff{\mbox{\sl Cliff}}
\def\ric{\mathop{\rm Ric}\nolimits}

\def\mod{\mathop{\rm mod}\nolimits}

\newcommand{\norm}[1]{\parallel\!\!{#1}\!\!\parallel}
\newcommand{\pin}{\hat{\circ}}
\newcommand{\htimes}{\hat{\otimes}}
\newcommand{\R}{\mathbb{R}}
\newcommand{\C}{\mathbb{C}}

\newcommand{\Z}{\mathbb{Z}}

\newcommand{\mf}{\mathfrak}
\newcommand{\mc}{\mathcal}
\newcommand{\dstar}{d\!\star\!}

\begin{document}

\title{\bf{Generalised $G_2$--manifolds}}

\author{Frederik Witt\\Freie Universit\"at Berlin\\FB Mathematik und Informatik\\Arnimallee 3\\D--14195 Berlin\\{\sf fwitt@math.fu--berlin.de}}

\date{}

\maketitle

\begin{abstract}
\noindent We define new Riemannian structures on 7--manifolds by a
differential form of mixed degree which is the critical point of a
(possibly constrained) variational problem over a fixed cohomology
class. The unconstrained critical points generalise the notion of a
manifold of holonomy $G_2$, while the constrained ones give rise to
a new geometry without a classical counterpart. We characterise
these structures by means of spinors and show the integrability
conditions to be equivalent to the supersymmetry equations on
spinors in type II supergravity theory with bosonic background
fields. In particular, this geometry can be described by two linear
metric connections with skew--symmetric torsion. Finally, we
construct explicit examples by introducing the device of $T$--duality.
\end{abstract}

\section{Introduction}

Over a 7--manifold, a topological reduction to a principal
$G_2$--bundle is achieved by the existence of a certain 3--form. The
fact that this 3--form is generic or {\em stable} (following the
language of Hitchin) enables one to set up a variational principle
over a fixed cohomology class whose critical points are precisely
the manifolds of holonomy $G_2$~\cite{hi01}.

In this paper we are concerned with a new type of Riemannian
geometry over 7--manifolds which generalises this notion.
Topologically speaking, it is defined by an even or odd form which
we think of as a {\em spinor} for the orthogonal bundle $T\oplus
T^*$ with its natural inner product of split signature. This
construction is perfectly general and works in all dimensions, but
there are special cases where $\R^*\times Spin(n,n)$ acts with an
{\em open} orbit in its spin representations $\Lambda^{ev,od}T^*$.
An example in dimension 6 are the so--called generalised Calabi--Yau
manifolds associated with an $SU(3,3)$--invariant spinor~\cite{hi03}.
In dimension 7, there are stable spinors whose stabiliser is either
conjugate to $G_2(\C)$ or $G_2\times G_2$. In this paper, we shall
deal with the latter case and define a generalised $G_2$--structure
as a (topological) reduction from $\R^*\times Spin(7,7)$ to $G_2\times G_2$.
Stability allows us to consider a generalised variational problem
which provides us with various integrability conditions.

We begin by introducing the algebraic setup. A reduction to
$G_2\times G_2$ gives rise to various objects. Firstly, it induces a
metric on $T$, and moreover, a 2--form $b$ which we refer to as the
{\em $B$--field}. As an element of the Lie algebra $\mf{so}(7)$
inside $\mf{so}(7,7)$, it acts on any $Spin(7,7)$--representation by exponentiation.
Secondly, we obtain two unit spinors $\Psi_+$ and $\Psi_-$ in the
irreducible spin representation of $Spin(7)$. Tensoring the spinors
yields an even or odd form $[\Psi_+\otimes\Psi_-]^{ev,od}$ by
projection on $\Lambda^{ev,od}$. The first important result we prove
states that any $G_2\times G_2$--invariant spinor $\rho$ can be
expressed as
\begin{equation}\label{genformrep}
\rho=e^{-\phi}\exp(b/2)\wedge[\Psi_+\otimes\Psi_-]^{ev,od}.
\end{equation}
In physicists' terminology, the scalar $\phi$ represents the {\em
dilaton} -- here it appears as a scaling factor.

Moving on to global issues, we see that up to a $B$--field
transformation, a generalised $G_2$--structure is essentially a pair
of principle $G_2$--fibre bundles inside the orthonormal frame bundle
determined by the metric. In particular, any topological
$G_2$--manifold trivially induces a generalised $G_2$--manifold, and
consequently, any spinnable seven--fold carries such a structure.
Over a compact manifold, we can classify generalised
$G_2$--structures up to vertical homotopy by an integer which
effectively counts (with an appropriate sign convention) the number
of points where the two $G_2$--structures inside the
$Spin(7)$--principal bundle coincide.

Over a closed manifold, we can then set up a variational problem
along the lines of~\cite{hi03}. In close analogy to the classical
case, the condition for a critical point is that both the even and
the odd $G_2\times G_2$--invariant spinor, now regarded as a form,
have to be closed. We shall adopt this as the condition of {\em
strong} integrability for any (not necessarily compact) generalised
$G_2$--structure. Interpreting this integrability condition in terms
of the right--hand side of~(\ref{genformrep}) leads to our main
result. Theorem~\ref{integrability} characterises strongly integrable
generalised $G_2$--structures in terms of two linear metric
connections $\nabla^{\pm}$ with skew--symmetric, closed torsion $\pm T$ such that
$$
(d\phi\pm\frac{1}{2}T)\cdot\Psi_{\pm}=0
$$
holds. Connections with skew--symmetric torsion have gained a lot of attention in the recent mathematical literature (see, for
instance,~\cite{agfr04},~\cite{friv02},~\cite{friv03} and
~\cite{iv01}) due to their importance in string theory.
Eventually, our reformulation yields the supersymmetry equations arising in type II supergravity with bosonic background fields~\cite{gmpw04}.

The spinorial picture is also useful for deriving geometrical properties.
In particular, we compute the Ricci tensor which is given by
$$
\ric(X,Y)=-2H^{\phi}(X,Y)+\frac{1}{4}g(X\llcorner
T,Y\llcorner T),
$$
with $H^{\phi}$ denoting the Hessian of the dilaton. A further striking
consequence is a no--go theorem for compact manifolds (generalising
similar statements in~\cite{ivpa01} and~\cite{gmpw04}): Here, the torsion of any strongly integrable generalised $G_2$--structure must vanish, that is, the underlying spinors $\Psi_+$ and $\Psi_-$ are
parallel with respect to the Levi--Civita connection. In this sense,
only ``classical" solutions can be found for the variational problem.
However, local examples with non--vanishing torsion exist in
abundance. Using again the form definition of a generalised
$G_2\times G_2$--structure, we will describe a systematic
construction method known in string theory as $T$--duality. It
consists of changing the topology by replacing a fibre isomorphic to
$S^1$ (or more generally to an $n$--torus) without destroying
integrability. This allows us to pass from a trivial generalised
structure coming from a classical $S^1$--invariant $G_2$--structure
with vanishing torsion (for which we can easily find examples) to a
trivial $S^1$--bundle with a non--trivial generalised $G_2$--structure.

The lack of interesting compact examples motivated us to consider a
{\em constrained} variational problem following ideas
in~\cite{hi01}. This gives rise to {\em weakly} integrable
structures of either {\em even} or {\em odd} type. For these the
no--go theorem does not apply, but the construction of examples, let
alone compact ones, remains an open problem. Although arising out of
a similar constraint as manifolds of weak holonomy $G_2$, this
notion gives rise to a new geometry without a classical counterpart
which renders a straightforward application of $T$--duality
impossible. We investigate its properties along with those of the
strongly integrable case.

The algebraic theory outlined above not only makes sense for $G_2\times
G_2$--structures, but also for $Spin(7)\times Spin(7)$ inside
$Spin(8,8)$, leading to the notion of a {\em generalised
$Spin(7)$--manifold}. It is also defined by an invariant form of
mixed degree, which, however, is not stable. Comparison with
classical $Spin(7)$--geometry suggests closeness of this form as a
natural notion of integrability. Section~\ref{spin7theory} briefly
explores the theory of these structures which is developed in full
detail in~\cite{wi04t}.

\medskip

This paper is based on a part of the author's doctoral
thesis~\cite{wi04t}. He wishes to acknowledge the DAAD, the Studienstiftung des deutschen Volkes, the EPSRC and the University of Oxford for various financial support. The author also wants to thank Gil Cavalcanti, Michael Crabb, Claus Jeschek, Wilson Sutherland and
the examiners Dominic Joyce and Simon Salamon for helpful comments
and discussions. Finally, he wishes to express his special gratitude to his
supervisor Nigel Hitchin.

\section{The linear algebra of generalised $G_2$--structures}

\subsection{Generalised metrics}\label{generalisedmetric}

We consider the vector bundle $T\oplus T^*$, where $T$ is a real,
$n$--dimensional vector space. It carries a natural orientation and
the inner product of signature $(n,n)$, defined for $v\in T$ and
$\xi\in T^*$ by
\begin{equation}\label{innprod}
(v+\xi,v+\xi)=-\frac{1}{2}\xi(v),
\end{equation}
singles out a group conjugate to $SO(n,n)$ inside $GL(2n)$. Note
that $GL(n)\leqslant SO(n,n)$. As a $GL(n)$--space, the Lie algebra
of $SO(n,n)$ decomposes as
$$
\mf{so}(n,n)=\Lambda^2(T\oplus
T^*)=End(T)\oplus\Lambda^2T^*\oplus\Lambda^2T.
$$
In particular, any 2--form $b$ defines an element in the Lie algebra
$\mf{so}(n,n)$. We will refer to such a 2--form as a {\em $B$--field}.
Exponentiated to $SO(n,n)$, its action on $T\oplus T^*$ is
$$
\exp(b)(v\oplus\xi)=v\oplus (v\llcorner b +\xi).
$$
Next we define an action of $T\oplus T^*$ on $\Lambda^*T^*$ by
$$
(v+\xi)\bullet\tau=v\llcorner\tau+\xi\wedge\tau.
$$
As this squares to minus the identity it gives rise to an
isomorphism $\cliff(T\oplus T^*)\cong End(\Lambda^*T^*)$. The
exterior algebra $S=\Lambda^*T^*$ becomes thus the pinor
representation space of $\cliff(T\oplus T^*)$ and splits into the
irreducible spin representation spaces $S^{\pm}=\Lambda^{ev,od}T^*$
of $Spin(n,n)$.

\begin{rmk}
There is a canonical embedding $GL_+(n)\hookrightarrow Spin(n,n)$ of the
identity component of $GL(n)$ into the spin group of $T\oplus T^*$.
As a $GL_+(n)$--module we have
$$
S^{\pm}=\Lambda^{ev,od}T^*\otimes(\Lambda^nT)^{1/2},
$$
in analogy to the complex case. There $U(n)\hookrightarrow
Spin^{\C}(2n)=Spin(2n)\times_{\Z_2}S^1$, and the even and odd forms
get twisted with the square root of the canonical line bundle. As
long as we are doing linear algebra this is a mere notational issue
but in the global situation we cannot trivialise $\Lambda^nT$ unless
the manifold is orientable. In fact, a more refined analysis reveals that we can always choose a spin structure for $T\oplus T^*$ -- whether the manifold is orientable or not -- such that the spinor bundle is isomorphic, as a $GL(n)$--space, to the exterior form bundle, albeit in a non--canonical way (see, for instance, the discussion in Section 2.8 in~\cite{gu03}). We will neglect these subtleties for we will consider orientable manifolds only, and therefore omit the twist to ease notation.
However it is important to bear it in mind when we set up the
variational formalism in Section~\ref{variation}.
\end{rmk}

Let $\sigma$ be the Clifford algebra anti--automorphism defined on
any element of degree $p$ by $\sigma(\alpha^p)=\epsilon(p)\alpha^p$
where $\epsilon(p)=1$ for $p\equiv0,3\mod4$ and $-1$ for
$p\equiv1,2\mod4$. The bilinear form
$$
\langle\alpha,\beta\rangle=\big(\alpha\wedge\sigma(\beta)\big)_n,
$$
where the subscript $n$ indicates taking the top degree component,
is non--degenerate and invariant under the action of $Spin(n,n)$. It
is symmetric if $n\equiv 0,3\mod4$ and skew if $n\equiv 1,2\mod4$.
Moreover, $S^+$ and $S^-$ are non--degenerate and orthogonal if $n$
is even and totally isotropic if $n$ is odd. Finally, we note that
the action of a $B$--field $b$ on a spinor $\tau$ is given by
$$
\exp(b)\bullet\tau=(1+b+\frac{1}{2}b\wedge
b+\ldots)\wedge\tau=e^{b}\wedge\tau.
$$
In this paper we shall be concerned with special structures on
$T\oplus T^*$, which we describe in terms of reductions to special
subgroups of $O(n,n)$ or $\R^*\times O(n,n)$.

\begin{definition}
A {\em generalised metric structure} is a reduction from $O(n,n)$ to
$O(n)\times O(n)$.
\end{definition}

\begin{figure}[ht]
\begin{center}
\input{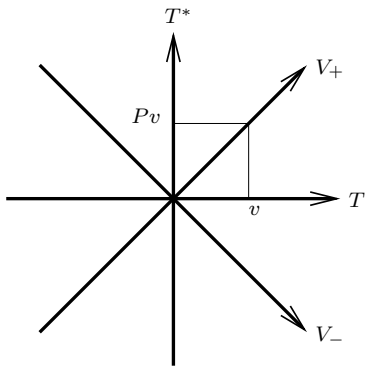}
\end{center}
\caption{Metric splitting of $T\oplus T^*$}\label{lightcone}
\end{figure}

Figure~\ref{lightcone} suggests how to characterise a metric
splitting algebraically.
If we think of the coordinate axes $T$ and $T^*$ as a lightcone,
choosing a subgroup conjugate to $O(n)\times O(n)$ inside $O(n,n)$
boils down to the choice of a spacelike $V_+$ and a timelike
orthogonal complement $V_-$. Interpreting $V_+$ as the graph of a
linear map $P_+:T\to T^*$ yields a metric $g$ and a 2--form $b$ as
the symmetric and the skew part of the corresponding bilinear form
$P_+\in T^*\otimes T^*$. Indeed we have
$$
g(t,t)=(t,P_+t)=(t\oplus P_+t,t\oplus P_+t)/2>0
$$
so that $g$ is positive definite. As $V_+$ and $V_-$ are orthogonal,
taking $V_-$ instead of $V_+$ yields the same 2--form $b$ but the
metric $-g$ . Conversely, assume we are given a metric $g$ and a
2--form $b$ on $T$. If we transform the diagonal $D_{\pm}=\{t\oplus
\mp t\llcorner g\:|\: t\in T \}$ by $\exp(b)$, we obtain a splitting
$V_+\oplus V_-$ inducing $g$ and $b$.

\begin{prp}\label{genmetricchar}
The choice of an equivalence class in the space $O(n,n)/O(n)\times
O(n)$ is equivalent to either set of the following data:

{\rm (i)} a metric splitting
$$
T\oplus T^*=V_+\oplus V_-
$$
into subbundles $(V_+,g_+)$ and $(V_-,g_-)$ with positive and
negative definite metrics $g_{\pm}=(\cdot\,,\cdot)_{|V_{\pm}}$.

{\rm (ii)} a Riemannian metric $g$ and a 2--form $b$ on $T$.
\end{prp}

Note that a reduction to $O(n)\times O(n)$ determines $g$ and $b$ up to a common scalar which, however, is fixed by the explicit choice of the inner product~(\ref{innprod}).

\begin{cor}
$$
O(n,n)/O(n)\times O(n)=\{P:T\to T^*\:|\:(Pt,t)>0\mbox{
for all }t\not= 0\}.
$$
\end{cor}

In the same vein, we call an element of $SO(n,n)/SO(n)\times SO(n)$
a {\em generalised oriented metric structure}, which corresponds to a
metric $g$, a $B$--field $b$ and an orientation on $T$. Since the bundle $T\oplus T^*$ is always spinnable, we can also lift the discussion to the group $Spin(n,n)$. Moreover, in some situations it is natural to introduce an additional scalar, that is, we enhance the structure group $Spin(n,n)$ to the conformal spin group $\R^*\times Spin(n,n)$ so that a reduction to $Spin(n)\times Spin(n)$ gives a further degree of freedom. This provides the right framework within which we can discuss generalised $G_2$--structures, to which we turn next.

\subsection{Generalised $G_2$--structures and stability}\label{geng2stru}

\begin{definition}
A {\em generalised $G_2$--structure} is a reduction from the
structure group $\R^*\times Spin(7,7)$ of $T\oplus T^*$ to $G_2\times G_2$.
\end{definition}

We want to characterise generalised $G_2$--structures along the lines
of Proposition~\ref{genmetricchar}. To get things rolling, we first
look at the tensorial invariants on $T$ which are induced by such a
reduction. Since $G_2\times G_2$ determines some group
$SO(V_+)\times SO(V_-)$ conjugate to $SO(7)\times SO(7)$, it induces
a generalised oriented metric structure $(g,b)$. The isomorphisms
$$
\pi_{b\pm}:x\in (T,\pm g)\mapsto -x\oplus
-x\llcorner(g\pm b)
$$
allow us to transport the respective $G_2$--structure on $V_+$ and $V_-$ to the tangent bundle. The resulting $G_2$--structures $G_{2+}$ and $G_{2-}$ inside $SO(T,g)=SO(7)$ give rise to two unit spinors $\Psi_+$ and $\Psi_-$ in the irreducible spin representation $\Delta=\R^8$ of
$Spin(7)=Spin(T,g)$.

On the other hand, $G_2\times G_2$ also acts on $\Lambda^*T^*$ as a
subgroup of $Spin(7,7)$. To relate these two actions we consider the
following construction which is basically the classical
identification of $\Delta\otimes\Delta$ with $\Lambda^*$ followed by
a twist with $\exp(b/2)$. We write the Clifford algebra
$\cliff(T\oplus T^*)$ as a $\Z_2$--graded tensor product
$\cliff(V_+)\htimes\cliff(V_-)\cong\cliff(T\oplus T^*)$ where the
isomorphism is given by extension of the map $v_+\htimes v_-\mapsto
v_+\bullet v_-$. The maps $\pi_{\pm}$ induce isomorphisms
$\cliff(T,\pm g)\cong\cliff(V_{\pm},g_{\pm})$ mapping
$Spin(7)=Spin(T,\pm g)$ isomorphically onto $Spin(V_{\pm},g_{\pm})$.
Consequently, the compounded algebra isomorphism
$$
\begin{array}{ccl}
\iota_b:\cliff(T,g)\htimes\cliff(T,-g) & \to &
\cliff(T\oplus T^*)\\
x\htimes y & \mapsto & \pi_{b+}(x)\bullet\pi_{b-}(y)
\end{array}
$$
maps $Spin(7)\times Spin(7)$ onto $Spin(V_+)\times Spin(V_-)$ inside
$Spin(7,7)$.

Let $q$ denote the $Spin(7)$--invariant inner product on
$\Delta$. For two spinors $\Psi_+$ and $\Psi_-$ the {\em pinor
product} $\Psi_+\pin\Psi_-$ is the endomorphism of $\Delta$ defined
by
$$
\Psi_+\pin\Psi_-(\Phi)=q(\Psi_-,\Phi)\Psi_+.
$$
On the other hand, there exists an algebra isomorphism $\kappa:\cliff(T)\to\End(\Delta)\oplus\End(\Delta)$. Projection on the first summand induces a matrix representation for Clifford multiplication. For concrete computations, we realise this representation by 
\begin{eqnarray}
e_1 & \mapsto & E_{1,2}-E_{3,4}-E_{5,6}+E_{7,8},\nonumber\\
e_2 & \mapsto & E_{1,3}+E_{2,4}-E_{5,7}-E_{6,8},\nonumber\\
e_3 & \mapsto & E_{1,4}-E_{2,3}-E_{5,8}+E_{6,7},\nonumber\\
e_4 & \mapsto & E_{1,5}+E_{2,6}+E_{3,7}+E_{4,8},\label{matrep}\\
e_5 & \mapsto & E_{1,6}-E_{2,5}+E_{3,8}-E_{4,7},\nonumber\\
e_6 & \mapsto & E_{1,7}-E_{2,8}-E_{3,5}+E_{4,6},\nonumber\\
e_7 & \mapsto & E_{1,8}+E_{2,7}-E_{3,6}-E_{4,5},\nonumber
\end{eqnarray}
where $E_{ij}=(\delta_{jk}\delta_{il}-\delta_{ik}\delta_{jl})_{k,l=1}^8$, $i<j$ is the standard basis of skew--symmetric endomorphisms of $\Delta$, taken with respect to an orthonormal basis $\Psi_1,\ldots,\Psi_8$. This relates to the pinor
product by
\begin{equation}\label{cliffordid}
(x\cdot\Psi_+)\pin\Psi_-=\kappa(x)\circ(\Psi_+\pin\Psi_-)\mbox{ and
}\Psi_+\pin(x\cdot\Psi_-)=(\Psi_+\pin\Psi_-)\circ\kappa\big(\sigma(x)\big).
\end{equation}
Let $j$ denote the canonical vector space isomorphism between
$\cliff(T)$ and $\Lambda^*T^*$. We consider the map
$$
[\cdot\,,\cdot\,]:\Delta\otimes\Delta\stackrel{\pin\oplus0}{\longrightarrow}
\End(\Delta)\oplus\End(\Delta)\stackrel{\kappa^{-1}}{\longrightarrow}\cliff(T,g)\stackrel{j}{\longrightarrow}
\Lambda^*T^*
$$
and think of any element in $\Delta\otimes\Delta$ as a form. We
denote by $[\Psi_+\otimes\Psi_-]^{ev,od}$ the projection on the even or
odd part and add a subscript $b$ if we wedge with the exponential $\exp(b/2)$. The following result states that up to a sign
twist, the action of $T$ on $\Delta\otimes\Delta$ and
$\Lambda^{ev,od}$ commute.

\begin{prp}\label{lmap}
For any $x\in T$ and $\Psi_+,\,\Psi_-\in\Delta$ we have
\begin{eqnarray}
\,[x\cdot\Psi_+\otimes\Psi_-]_b^{ev,od} & = &
\phantom{\pm}\iota_b(x\htimes1)\bullet \,[\Psi_+\otimes\Psi_-]_b^{od,ev}\nonumber\\
\,[\Psi_+\otimes x\cdot\Psi_-]_b^{ev,od} & = & \pm\iota_b(1\htimes
x)\bullet [\Psi_+\otimes\Psi_-]_b^{od,ev}.\nonumber
\end{eqnarray}
In particular, the forms $[\Psi_+\otimes\Psi_-]^{ev,od}_b$ are $G_2\times G_2$--invariant.
\end{prp}

\begin{prf}
First we assume that $b\equiv 0$. By convention, we let $T$ act
through the inclusion
$$
T\hookrightarrow\cliff(T)\stackrel{
\kappa}{\cong}\End(\Delta)\oplus\End(\Delta)
$$
followed by projection on the first summand. Thus
\begin{eqnarray*}
\,[x\cdot\Psi_+\otimes\Psi_-] & = &
j\big(\kappa^{-1}(x\cdot\Psi_+\pin\Psi_-\oplus0)\big)\\
& = & j\big(x\cdot \kappa^{-1}(\Psi_+\pin\Psi_-)\big)\\
& = & -x\llcorner [\Psi_+\otimes\Psi_-]+x\wedge
[\Psi_+\otimes\Psi_-]\\
& = & \iota_0(x\htimes 1)\bullet [\Psi_+\otimes\Psi_-],
\end{eqnarray*}
where we have used~(\ref{cliffordid}) and the identity $j(v\cdot
x)=v\wedge x -v\llcorner x$ which holds for any Clifford algebra
$\cliff(V)$ and $x\in V$. We can argue similarly for $\Psi_+\otimes
x\cdot\Psi_-$ and we obtain
$$
[x\cdot\Psi_+\otimes\Psi_-]^{ev}=\iota_0(x\htimes 1)\bullet
[\Psi_+\otimes\Psi_-]^{od},\quad
[x\cdot\Psi_+\otimes\Psi_-]^{od}=\iota_0(x\htimes 1)\bullet
[\Psi_+\otimes\Psi_-]^{ev},
$$
so that the case $b=0$ is shown.

Now let $b$ be an arbitrary $B$--field. For the sake of clarity we
will temporarily denote by $\widetilde{\exp}$ the exponential map
from $\mf{so}(7,7)$ to $Spin(7,7)$, while the untilded exponential takes
values in $SO(7,7)$. The adjoint representation $Ad$ of the group of
units inside a Clifford algebra $\cliff(V)$ restricts to the double
cover $Spin(V)\to SO(V)$ still denoted by $Ad$. As a transformation
in $SO(V)$ we then have $Ad\circ\widetilde{\exp}=\exp\circ\, ad$.
Since $ad([v,w])=4v\wedge w$ and the vectors $e_1,\ldots, e_7$ are isotropic with respect to the inner product $\langle\cdot\,,\cdot\,\rangle$, in our situation we get $ad(e_i\bullet e_j)=2e_i\wedge e_j$ and thus
$$
e^b=Ad(\widetilde{e}^{\sum b_{ij}e_i\bullet e_j/2})
$$
for the $B$--field $b=\sum_{i<j}b_{ij}e_i\wedge e_j$. Hence
\begin{eqnarray*}
[x\cdot\Psi_+\otimes\Psi_-]^{ev,od}_b & = &
\big(\widetilde{e}^{\:b/2}\bullet\iota_0(x\htimes
1\big)\bullet\widetilde{e}^{\:-b/2})\bullet [\Psi_+\otimes\Psi_-]^{od,ev}_b\\
& = & Ad(\widetilde{e}^{b/2})\big(\iota_0(x\htimes 1)\big)\bullet
[\Psi_+\otimes\Psi_-]^{od,ev}_b\\
& = & \iota_b(x\htimes 1)\bullet [\Psi_+\otimes\Psi_-]^{od,ev}_b.
\end{eqnarray*}
Similarly, the claim is checked for $[\Psi_+\otimes
x\cdot\Psi_-]^{ev,od}_b$ which completes the proof.
\end{prf}

The identification $\Delta\otimes\Delta\cong\Lambda^{ev,od}$ also
enables us to derive a normal form for
$[\Psi_+\otimes\Psi_-]^{ev,od}$ by choosing a suitable orthonormal
frame. The coefficients of the homogeneous components are given up
to a scalar by $q(\kappa(e_I)\cdot\Psi_+,\Psi_-)$. Since
the action of $Spin(7)$ on the Stiefel variety $V_2(\Delta)$, the
set of pairs of orthonormal spinors, is transitive, we may assume
that $\Psi_+=\Psi_1$ and $\Psi_-=c\Psi_1+s\Psi_2$, where $c$ and $s$
is shorthand for $\cos(a)$ and $\sin(a)$ with
$a=\sphericalangle(\Psi_+,\Psi_-)$. Using the
representation~(\ref{matrep}), the computation of the normal form is
straightforward.

\begin{prp}\label{normalform}
There exists an orthonormal basis $e_1,\ldots,e_7$ such that
\begin{eqnarray*}
[\Psi_+\otimes\Psi_-]^{ev} & = &
\phantom{+}\cos(a)+\sin(a)(-e_{23}-e_{45}+e_{67})+\\
& &
+\cos(a)(-e_{1247}+e_{1256}+e_{1346}+e_{1357}-e_{2345}+e_{2367}+e_{4567})+\\
& & +\sin(a)(e_{1246}+e_{1257}+e_{1347}-e_{1356})-\sin(a)e_{234567}
\end{eqnarray*}
and
\begin{eqnarray*}
[\Psi_+\otimes\Psi_-]^{od} & = &
\phantom{+}\sin(a)e_1+\sin(a)(e_{247}-e_{256}-e_{346}-e_{357})+\\
& &
+\cos(a)(e_{123}+e_{145}-e_{167}+e_{246}+e_{257}+e_{347}-e_{356})+\\
& &
+\sin(a)(-e_{12345}+e_{12367}+e_{14567})+\cos(\alpha)e_{1234567}.
\end{eqnarray*}
\end{prp}

If the spinors $\Psi_+$ and $\Psi_-$ are linearly independent, we
can express the structure form in terms of the invariants associated
with $SU(3)$, the stabiliser of the pair $(\Psi_1,\Psi_2)$.

\begin{cor}\label{normalformcor}
If $\alpha$ denotes the dual of the unit vector in $T$ which is
stabilised by $SU(3)$, $\omega$ the K\"ahler form and $\psi_{\pm}$
the real and imaginary parts of the holomorphic volume form, then
$$
[\Psi_+\otimes\Psi_-]^{ev}=c+s\omega+c(\alpha\wedge\psi_--\frac{1}{2}\omega^2)-s\alpha\wedge\psi_+-\frac{1}{6}s\omega^3
$$
and
$$
[\Psi_+\otimes\Psi_-]^{od}=s\alpha-c(\psi_++\alpha\wedge\omega)-s\psi_--\frac{1}{2}s\alpha\wedge\omega^2+cvol_g.
$$
Moreover, $\Psi_-=s\alpha\cdot\Psi_+$. If $a=0$, then both
$G_2$--structures coincide and
$$
[\Psi_+\otimes\Psi_-]^{ev}=1-\star\varphi,\quad[\Psi_+\otimes\Psi_-]^{od}=-\varphi+vol_g,
$$
where $\varphi$ is the stable 3--form associated with $G_2$.
\end{cor}

\begin{rmk}
The stabiliser of the forms $1+\star\varphi$ and $\varphi+vol_g$ is
isomorphic to $G_2(\C)$~\cite{wi04t}.
\end{rmk}

The $63$ degrees of freedom which parametrise a reduction from $Spin(7,7)$ to $G_2\times G_2$ are exhausted by $28$ degrees of freedom for the choice of a metric $g$, $21$ for a 2--form $b$ and twice $7$ degrees for two unit spinors $\Psi_+$ and $\Psi_-$. However, this data does not achieve a full description of $G_2\times G_2$--invariant spinors yet as these are {\em stable}, a notion due to Hitchin~\cite{hi03}.

\begin{definition}
A spinor $\rho$ in $\Lambda^{ev,od}T^{n*}$ is said to be {\em
stable} if $\rho$ lies in an open orbit under the action of
$\R^*\times Spin(n,n)$.
\end{definition}

In~\cite{saki77} Sato and Kimura classified the representations of
complex reductive Lie groups which admit an open orbit. Apart from
dimension 7, stable spinors in the sense of the definition above
only occur in dimension 6 and give rise to {\em generalised
Calabi--Yau--structures} associated with the group
$SU(3,3)$~\cite{hi03}. The key point here is that in both cases we
obtain an invariant volume form. In our case, this generalises the
concept of stability for a $G_2$--invariant 3--form $\varphi$ with
associated volume form
$\varphi\wedge\star_{\varphi}\varphi$~\cite{hi01}.

To begin with, we note that the spin representation $\Delta$ of
$Spin(7)$ is real, and so is the tensor product
$\Delta\otimes\Delta$. Consequently, there is (up to a scalar) a
unique invariant in $\Lambda^{ev}\otimes\Lambda^{od}$, or
equivalently, $Spin(7)\times Spin(7)$--equivariant maps
$\Lambda^{ev,od}T^*\to\Lambda^{od,ev}T^*$. Morally these are given
by the Hodge $\star$--operator twisted with the $B$--field and the
anti--automorphism $\sigma$.

\begin{definition}
The {\em box--operator} or {\em generalised Hodge $\star$--operator}
$\Box_{g,b}:\Lambda^{ev,od}T^*\to\Lambda^{od,ev}T^*$ associated with the generalised metric $(g,b)$ is defined by
$$
\Box_{g,b}\rho=e^{b/2}\wedge\star_g\sigma(e^{-b/2}\wedge\rho).
$$
\end{definition}

If $g$ and $b$ are induced by $\rho$ we will also use the sloppier
notation $\Box_{\rho}$ or drop the subscript altogether. For present
and later use, we note the following lemma whose proof is immediate
from the definitions.

\begin{lem}\label{hat}
Let $\rho\in\Lambda^{ev,od}T^*$. Then
$\sigma(\star\rho)=\star\sigma(\rho)$ and
$\sigma(e^b\wedge\rho)=e^{-b}\wedge\sigma(\rho)$.
\end{lem}

\begin{prp}\label{selfdual}
For any $\Psi_+,\Psi_-\in\Delta\otimes\Delta$
$$
\Box_{g,b}[\Psi_+\otimes\Psi_-]_b=[\Psi_+\otimes\Psi_-]_b
$$
or equivalently,
$$
\Box_{g,b}[\Psi_+\otimes\Psi_-]_b^{ev,od}=[\Psi_+\otimes\Psi_-]_b^{od,ev}.
$$
\end{prp}

\begin{prf}
According to Lemma~\ref{hat},
\begin{eqnarray*}
\Box_{g,b}[\Psi_+\otimes\Psi_-]_b
& = & e^{b/2}\wedge\star\sigma[\Psi_+\otimes\Psi_-]\\
& = & e^{b/2}\wedge j\big(\kappa^{-1}(\Psi_+\pin\Psi_-)\cdot vol\big)
\end{eqnarray*}
where we used the general identity $\star j(x)= j\big(\sigma(x)\cdot
vol\big)$. But $\kappa^{-1}\big((\Psi_+\pin\Psi_-\oplus0)\cdot vol\big)$
is just $ \kappa^{-1}(\Psi_+\pin\Psi_-\oplus0)$, whence the assertion.
\end{prf}

Let $U$ denote the space of $G_2\times G_2$--invariant spinors in
$\Lambda^{ev}$ or $\Lambda^{od}$. With any element of $U$ we
associate a volume form
$$
\mc{Q}:\rho\in U\mapsto q(\Box_{\rho}\rho,\rho)\in\Lambda^7T^*,
$$
Now the $\Box$--operator transforms naturally under the lift
$\widetilde{A}\in Pin(7,7)$ of any element $A\in O(7,7)$ which means
that
$$
\Box_{\widetilde{A}\bullet\rho}\widetilde{A}\bullet\rho=\widetilde{A}\bullet\Box_{\rho}\rho.
$$
Therefore, we immediately conclude

\begin{prp}\label{volformg2g2}
$\mc{Q}$ is homogeneous of degree 2 and $Spin(7,7)$--invariant.
\end{prp}

\begin{rmk}
An explicit coordinate description of the complexified invariant was
given in~\cite{gy90}. However, with the aim of setting up the
variational principle, this formulation proved to be rather
cumbersome for our purposes, which motivated our approach in terms of
$G$--structures.
\end{rmk}

We will also need the differential of this map. Since the form $\langle\cdot\,,\cdot\,\rangle$ is non--degenerate, we can write
$$
D\mc{Q}_{\rho}(\dot{\rho})=\langle\hat{\rho},\dot{\rho}\rangle.
$$
for a unique $\hat{\rho}\in\Lambda^{od}T^*$, the {\em
companion} of $\rho$, which is also a $G_2\times G_2$--invariant
spinor. By rescaling $\mc{Q}$ appropriately, we conclude that
$\hat{\rho}=\Box_{\rho}\rho$.

From $\mc{Q}$ we derive a further invariant attached to a $G_2\times
G_2$--invariant spinor $\rho$, namely a real scalar $\phi_{\rho}$
which we refer to as the {\em dilaton}. It is defined by
$$
\mc{Q}(\rho)=8e^{-2\phi_{\rho}}vol_g.
$$
Using the normal form description of Proposition~\ref{normalform} we
obtain $\mc{Q}\big([\Psi_+\otimes\Psi_-]^{ev,od}_b\big)=8vol_g$,
hence $\rho=e^{-\phi_{\rho}}[\Psi_+\otimes\Psi_-]^{ev,od}_b$.

We summarise the results of this section in the following

\begin{thm}\label{representation}
Generalised $G_2$--structures are in 1--1 correspondence with lines of
spinors $\rho$ in $\Lambda^{ev}T^*$ (or $\Lambda^{od}T^*$) whose
stabiliser under the action of $Spin(7,7)$ is isomorphic to
$G_2\times G_2$. We refer to $\rho$ as the {\em structure form} of
the generalised $G_2$--structure. This form can be uniquely written
(modulo a simultaneous sign change for $\Psi_+$ and $\Psi_-$) as
$$
\rho=e^{-\phi}[\Psi_+\otimes\Psi_-]^{ev}_b
$$
for a 2--form $b$, two unit spinors $\Psi_{\pm}\in\Delta$ and a real
scalar $\phi$. In particular, the space of $G_2\times G_2$--invariant spinors is open.
\end{thm}

\section{Topological generalised $G_2$--structures}\label{topgeng2}

\begin{definition}
A {\em topological generalised $G_2$--structure} over a 7--manifold
$M$ is a topological $G_2\times G_2$--reduction of the
$\R^*\times Spin(7,7)$--principal bundle associated with $T\oplus T^*$. It is
characterised by a stable even or odd spinor $\rho$ which we view as
a form. Consequently, we will denote this structure by the pair
$(M,\rho)$ and call $\rho$ the {\em structure form}.
\end{definition}

We will usually drop the adjective ``topological" and simply refer
to a generalised $G_2$--structure. As we saw earlier, a generalised
$G_2$--structure induces a generalised metric. In particular there
exists a metric $g$ or equivalently an $SO(7)$--principal fibre
bundle which admits two $G_2$--subbundles. The inclusion
$G_{2\pm}\subset Spin(7)$ implies that the underlying manifold is
spinnable and distinguishes a preferred spin structure for which we
consider the associated spinor bundle $\Delta$. Using
Theorem~\ref{representation} we can now assert the following
statement.

\begin{thm}\label{topology}
A topological generalised $G_2$--structure $(M,\rho)$ is
characterised by the following data:
\begin{itemize}
\item an orientation
\item a metric $g$
\item a 2--form $b$
\item a scalar function $\phi$
\item two unit spinors $\Psi_+,\Psi_-\in\Delta$ such that
$$
e^{-\phi}[\Psi_+\otimes\Psi_-]_b=\rho+\Box_{g,b}\rho.
$$
\end{itemize}
\end{thm}

A trivial example of a generalised $G_2$--structure is a topological
$G_2$--manifold with associated spinor $\Psi=\Psi_{\pm}$. A
generalised structure arising in this way (possibly with a non--trivial
$B$--field and dilaton) is said to be {\em straight}. The existence of a nowhere
vanishing spinor field in dimension 7 is a classical
result~\cite{lami89} and implies

\begin{cor}
A 7--fold $M$ carries a topological generalised $G_2$--structure if
and only if $M$ is spinnable.
\end{cor}

Next we discard the $B$--field and the dilaton and focus on the
$G_2$--structures induced by $\Psi_+$ and $\Psi_-$. Our aim is to
classify the $G_2$--structures up to equivalence under a
$Spin(7)$--fibre bundle isomorphism. Since the classification of
principal fibre bundles is a problem of homotopy theory, this boils
down to deform homotopically $\Psi_+$ into $\Psi_-$ through
sections. More concretely, we regard $G_2$--structures as being
defined by (continuous) sections of the sphere bundle
$p_{\mathbb{S}}:\mathbb{S}\to M$ associated with $\Delta$. On the space of
sections $\Gamma(\mathbb{S})$ we introduce the following equivalence
relation. Two spinors $\Psi_+$ and $\Psi_-$ are
considered to be equivalent (denoted $\Psi_+\sim\Psi_-$) if and only
if there exists a continuous map $G:M\times I\to \mathbb{S}$ such
that $G(x,0)=\Psi_+(x)$, $G(x,1)=\Psi_-(x)$ and $p_{\mathbb{S}}\circ
G(x,t)=x$. An equivalence class will be denoted by $[\Psi]$. If two
sections are vertically homotopic, then the corresponding
$G_2$--structures are isomorphic as principal $G_2$--bundles over $M$.
In particular the generalised structure defined by the pair
$(\Psi_+,\Psi_-)$ is equivalent to a straight structure if and only
if $\Psi_+\sim\Psi_-$. What we aim to determine is the set of
generalised structures with fixed $\Psi_+$, i.e.
$Gen(M)=\Gamma(\mathbb{S})/\sim$, the set of isomorphism classes of
principal $G_2$--fibre bundles. If a generalised structure is defined
by two inequivalent spinors, then it is said to be {\em exotic}. Here is an example.

\begin{ex}
Consider $M=S^7$. Since the tangent bundle of $S^7$ is trivial, so
is the sphere bundle of $\Delta$, i.e. $\mathbb{S}=S^7\times S^7$.
Consequently,
$$
Gen(S^7)=\Gamma(\mathbb{S})/\!\sim\;=[S^7,S^7]=\pi_7(S^7)=\Z
$$
and any map $\Psi_-:S^7\to S^7$ which is not homotopic to a constant
gives rise to an exotic structure for
$\Psi_+\equiv const$.
\end{ex}

In general, the question whether or not two sections are
vertically homotopic can be tackled by using obstruction theory (see
for instance~\cite{st51}). Assume that we are given a
fibre bundle with connected fibre over a not necessarily compact
$n$--fold $M^n$, and two sections $s_1$ and $s_2$ which are
vertically homotopy equivalent over the $q$--skeleton $M^{[q]}$ of
$M$. The obstruction for extending the vertical homotopy to the
$q+1$--skeleton lies in $H^{q+1}(M,\pi_{q+1}(F))$. In particular,
there is the first non--trivial obstruction $\delta(s_1,s_2)\in
H^m(M,\pi_m(F))$ for the least integer $m$ such that $\pi_m(F)\not=0$, called the {\em primary difference} of $s_1$ and $s_2$. It is a
homotopy invariant of $s_1$ and $s_2$ and enjoys the
additivity property
\begin{equation}\label{additivity}
\delta(s_1,s_2)+\delta(s_2,s_3)=\delta(s_1,s_3).
\end{equation}
Coming back to the generalised $G_2$--case we consider the sphere
bundle $\mathbb{S}$ over $M^7$ with fibre $S^7$. Consequently, the
primary difference of two sections lies in $H^7(M,\Z)$ and this is
the only obstruction for two sections to be vertically homotopy
equivalent. The additivity property implies that
$\delta(\Psi,\Psi_1)=\delta(\Psi,\Psi_2)$ if and only if
$\delta(\Psi_1,\Psi_2)=0$, that is, $\Psi_1\sim\Psi_2$. Moreover,
for any class $d\in H^7(M,\Z)$ there exists a section $\Psi_d$ such
that $d=\delta(\Psi,\Psi_d)$~\cite{st51}. As a consequence, we
obtain the

\begin{prp}
The set of generalised $G_2$--structures can be identified with
$$
Gen(M)=H^7(M,\Z)=\left\{\begin{array}{ll} \Z, & \mbox{ if $M$ is
compact}\\ 0, & \mbox{ if $M$ is non--compact}\end{array}\right..
$$
\end{prp}

Modulo $(b,\phi)$, generalised $G_2$--structures are therefore classified by an integer invariant which over a compact $M^7$ has the natural interpretation
as the number of points (counted with an appropriate sign
convention) where the two $G_2$--structures coincide. To see this we
associate with every equivalence class $[\Psi_-]$ the intersection
class $\#\big(\Psi_+(M),\Psi_-(M)\big)\in H_{14}(\mathbb{S},\Z)$ of
the 7--dimensional oriented submanifolds $\Psi_+(M)$ and $\Psi_-(M)$
inside $\mathbb{S}$. Since the total space of the sphere bundle is
14--dimensional, the intersection class counts the number of points in
$M$ where the two spinors $\Psi_+$ and $\Psi_-$ coincide. Taking the
cup product of the Poincar\'e duals of $\Psi_-(M)$ and $\Psi_+(M)$
sets up a map
$$
[\Psi_-]\in Gen(M)\mapsto PD\big(\Psi_+(M)\big)\cup
PD\big(\Psi_-(M)\big)\in H^{14}(\mathbb{S},\Z).
$$
On the other hand, integration along
the fibre defines an isomorphism
$\pi_{\mathbb{S}*}:H^{14}(\mathbb{S},\Z)\to H^7(M,\Z)$ by the Gysin sequence. Therefore, any generalised $G_2$--structure induced by the
equivalence class $[\Psi_-]$ over a compact 7--manifold $M$ gives rise to
a well--defined cohomology class $d(\Psi_+,\Psi_-)\in H^7(M,\Z)$. The
following theorem shows this class to coincide with the primary
difference. For the proof, I benefited from discussions with
W.~Sutherland and M.~Crabb.

\begin{thm}\label{intersection}
If $M$ is compact, then
$$
d(\Psi_+,\Psi_-)=\delta(\Psi_+,\Psi_-).
$$
In particular generalised $G_2$--structures are classified by the number
of points where the two underlying $G_2$--structures coincide.
\end{thm}

\begin{prf}
We regard the spinor bundle $\Delta$ as an 8--dimensional oriented
real vector bundle over M and consider the two sections $\Psi_+$ and
$\Psi_-$ of the sphere bundle. Deform the section $(x,t)\mapsto(1 -
t)\Psi_+(x) + t\Psi_-(x)$ of the pullback of $\Delta$ to $M\times\R$ to
be transverse to the zero--section. The primary difference
$\delta(\Psi_+,\Psi_-)$ can then be represented geometrically by the
(signed) zero--set of this deformation. In particular, if $\Psi_+$
and $-\Psi_-$ never coincide, then the primary difference is 0.
Therefore the intersection number, defined geometrically by making
$\Psi_+(M)$ and $\Psi_-(M)$ transverse and taking the coincidence
set, will be $\delta(\Psi_+,-\Psi_-)$ (with appropriate sign
conventions). By virtue of~(\ref{additivity}), we have
$\delta(\Psi_+,-\Psi_-) =
\delta(\Psi_+,\Psi_-)+\delta(\Psi_-,-\Psi_-)$. The difference class
$\delta(\Psi_-,-\Psi_-)$ corresponds to the self--intersection number
$\#(\Psi_-(M),\Psi_-(M))$ which is 0 since $M$ is 7--dimensional,
whence the assertion.
\end{prf}

\section{Integrable generalised $G_2$--structures and supersymmetry}\label{geometry}

\subsection{A variational principle}\label{variation}

Assume $M$ to be a closed and oriented manifold which carries a
topological generalised $G_2$--structure defined by the stable spinor
$\rho$. We shall set up a variational problem along the lines of
\cite{hi01} and \cite{hi03}. From a $GL_+(7)$--point of view, $\rho$
is a section of the vector bundle whose fibre is
$\Lambda^{ev,od}T^*\otimes(\Lambda^7T)^{1/2}$ (cf. the remark in
Section~\ref{generalisedmetric}). Untwisting by the line bundle
$(\Lambda^7T)^{1/2}$, we obtain a corresponding open set still
denoted by $U$, on which we can consider the induced volume
functional $\mc{Q}$ as described in Proposition~\ref{volformg2g2}.
It is homogeneous of degree 2 and therefore, it defines a
$GL(7)_+$--equivariant function
$\mc{Q}:U\subset\Lambda^{ev,od}T^*\to\Lambda^7T^*$ since
$$
\mc{Q}(A^*\rho)=(\det A)^{-1}\mc{Q}(A\bullet\rho)=(\det
A)^{-1}\mc{Q}(\rho).
$$
Associated with the $GL_+(7)$--principal fibre bundle over $M$,
$\mc{Q}$ thus takes values in $\Omega^7(M)$ and we obtain the volume
functional
$$
V(\rho)=\int_M\mc{Q}(\rho)
$$
defined over stable forms. As stability is an open condition, we can
differentiate this functional and consider its variation over a
fixed cohomology class. Instead of working with ordinary cohomology
only we will allow for an extra twist by a closed 3--form $H$. This
means that we replace the differential operator $d$ by the twisted
operator $d_H=d+H\wedge$. Closeness of $H$ guarantees that $d_H$
still defines a differential complex. Moreover, we will also
consider the following constraint. The bilinear form $\langle\cdot\,,\cdot\,\rangle$ induces a
non--degenerate pairing between $\Omega^{ev}(M)$
and $\Omega^{od}(M)$ defined by integration of $\langle\eta,\tau\rangle$ over
$M$. If $\eta=d_H\gamma$ is $H$--exact, Stokes' theorem and the
definition of the involution $\sigma$ imply that
\begin{equation}\label{stokes}
\int_M\langle d_H\gamma,\tau\rangle=\int_M\langle\gamma,d_H\tau\rangle.
\end{equation}
This vanishes for all $\gamma$ if and only if $\tau$ is $H$--closed.
Consequently, we can identify
$$
\Omega^{ev}_{\mbox{\tiny
$H$--exact}}(M)^*\cong\Omega^{od}(M)/\Omega^{od}_{\mbox{\tiny
$H$--closed}}(M).
$$
The exterior differential $d_H$ maps the latter space isomorphically
onto $\Omega^{ev}_{\mbox{\tiny $H$--exact}}(M)$ so that
$$
\Omega^{ev}_{\mbox{\tiny $H$--exact}}(M)^*\cong\Omega^{ev}_{\mbox{\tiny
$H$--exact}}(M).
$$
Finally, we obtain the non--degenerate quadratic form on
$\Omega^{ev}_{\mbox{\tiny $H$--exact}}(M)$ given by
$$
Q(d_H\gamma)=\int_M\langle\gamma,d_H\gamma\rangle.
$$
The same conclusion holds for odd instead of even forms.

\begin{thm}
Let $H$ be a closed 3--form.

{\rm (i)} A $d_H$--closed stable form
$\rho\in\Omega^{ev,od}(M)$ is a critical point in its cohomology
class if and only if $d_H\hat{\rho}=0$.

{\rm (ii)} A $d_H$--exact form $\rho\in\Omega^{ev,od}(M)$ defines a
critical point subject to the constraint $Q(\rho)=const$ if and
only if there exists a real constant $\lambda$ with
$d_H\rho=\lambda\hat{\rho}$.
\end{thm}

\begin{prf}
The first variation of $V$ is
$$
\delta
V_{\rho}(\dot{\rho})=\int_MD\mc{Q}_{\rho}(\dot{\rho})=\int_M\langle\hat{\rho},\dot{\rho}\rangle.
$$
To find the unconstrained critical points we have to vary over a
fixed $d_H$--cohomology class, i.e. $\dot{\rho}=d_H\gamma$. As we saw
in (\ref{stokes}), we have
$$
\delta
V_{\rho}(\dot{\rho})=\int_M\langle\hat{\rho},d_H\gamma\rangle=\int\limits_M\langle d_H\hat{\rho},\gamma\rangle,
$$
and this vanishes for all $\gamma$ if and only if $\hat{\rho}$ is
$d_H$--closed. In the constrained case the differential of $Q$
at $\rho$ is
$$
(\delta Q)_{\rho}(d_H\gamma)=2\int_M\langle\rho,\gamma\rangle.
$$
By Lagrange's theorem, we see that for a critical point we need
$d_H\hat{\rho}=\lambda\rho$.
\end{prf}

We adopt these various conditions for defining a critical point as
integrability condition of a topological generalised
$G_2$--structure.

\begin{definition}\label{intstruc}
Let $H$ be a closed 3--form and $\lambda$ be a real, non--zero
constant.

{\rm (i)} A topological generalised $G_2$--structure $(M,\rho)$ is
said to be {\em strongly integrable} with respect to $H$ if and only
if
$$
d_H\rho=0,\quad d_H\hat{\rho}=0.
$$

{\rm (ii)} A topological generalised $G_2$--structure $(M,\rho)$ is
said to be {\em weakly integrable with respect to $H$ and with
Killing number $\lambda$} if and only if
$$
d_H\rho=\lambda\hat{\rho}.
$$
We call such a structure {\em even} or {\em odd} according to the
parity of the form $\rho$. If we do not wish to distinguish the
type, we will refer to both structures as {\em weakly integrable}.

Similarly, the structures in {\rm (i)} and {\rm (ii)} will also be referred to
as {\em integrable} if a condition applies to both weakly and
strongly integrable structures.
\end{definition}

As we shall see in Corollary~\ref{torsioncomp}, the number $\lambda$
represents the $0$--torsion form of the underlying $G_2$--structures
which is why we refer to it as the Killing number~\cite{fkms97}.

\begin{ex}\hfill\newline
Consider a straight topological $G_2$--manifold $(M,\varphi,b,\phi)$
equipped with an additional closed 3--form $H$. According to Corollary~\ref{normalformcor} and Theorem~\ref{representation}, the corresponding structure form is
$\rho=e^{-\phi}e^{b/2}\wedge(1-\star\varphi)$ with companion
$\hat{\rho}=e^{-\phi}e^{b/2}\wedge(-\varphi+vol_g$). Writing
$T=db/2+H$, we want to solve the equations of strong integrability
$$
d_Te^{-\phi}(1-\star\varphi)=0,\quad d_Te^{-\phi}(-\varphi+vol_g)=0.
$$
It follows immediately that this is equivalent to $d\phi=0$, $T=0$
and $d\varphi=0$, $\dstar\varphi=0$, that is, the holonomy is
contained in $G_2$. If we ask for weak integrability the question
only makes sense in the even case as $\cos(a)=0$ for structures of
odd type. The equation of weak integrability becomes
$$
d_T(1-\star\varphi)-d\phi\wedge(1-\star\varphi)=\lambda(-\varphi+vol_g),
$$
implying $d\phi=0$, $T=-\lambda\varphi$ and
$-\lambda\varphi\wedge\star\varphi=\lambda vol_g$. Since
$\varphi\wedge\star\varphi=7vol_g$, we have $\lambda=0$. A straight
structure can therefore never induce a weakly integrable structure
-- in this sense, weak integrability has no classical counterpart.
\end{ex}

\subsection{Spinorial solution of the variational problem and supergravity}

Next we want to interpret the integrability conditions in
Definition~\ref{intstruc} in terms of the data of
Theorem~\ref{topology}. For a vanishing $B$--field, the
identification $\Delta\otimes\Delta\cong\Lambda^{ev,od}$ transforms
the twisted Dirac operators $\mc{D}$ and $\mc{\widehat{D}}$ on
$\Delta\otimes\Delta$, given locally by
\begin{eqnarray*}
\mc{D}(\Psi_-\otimes\Psi_+) 
& = & 
\sum e_i\cdot\nabla_{e_i}\Psi_-\otimes\Psi_++e_i\cdot\Psi_-\otimes\nabla_{e_i}\Psi_+,\\
\widehat{\mc{D}}(\Psi_-\otimes\Psi_+)
& = & 
\sum \nabla_{e_i}\Psi_-\otimes e_i\cdot\Psi_++\Psi_-\otimes e_i\cdot\nabla_{e_i}\Psi_+
\end{eqnarray*}
into the Dirac operators on $p$--forms $d+d^*$ and $(-1)^p(d\pm
d^*)$~\cite{lami89}. Here and in the sequel, $\nabla$ designates the
Levi--Civita connection on the tangent bundle $T$ or its lift to
$\Delta$. The transformation under a non--trivial $B$--field is given
in the next proposition.

\begin{prp}\label{dirac}
We have
\begin{eqnarray*}
\,[\mc{D}(\Psi_+\otimes\Psi_-)]^{ev,od}_b
& = &
\phantom{+}d[\Psi_+\otimes\Psi_-]^{od,ev}_b+d^{\Box}[\Psi_+\otimes\Psi_-]^{od,ev}_b\\
& &+\frac{1}{2}e^{b/2}\wedge\big(db\llcorner
[\Psi_+\otimes\Psi_-]^{od,ev}-db\wedge [\Psi_+\otimes\Psi_-]^{od,ev}\big)\\
\,[\widehat{\mc{D}}(\Psi_+\otimes\Psi_-)]^{ev,od}_b
& = &
\pm d[\Psi_+\otimes\Psi_-]^{od,ev}_b\mp d^{\Box}[\Psi_+\otimes\Psi_-]^{od,ev}_b\\
& & 
\mp\frac{1}{2}e^{b/2}\wedge\big(db\llcorner[\Psi_+\otimes\Psi_-]^{od,ev}+ db\wedge [\Psi_+\otimes\Psi_-]^{od,ev}\big),
\end{eqnarray*}
where $d^{\Box}\rho=\Box\, d\,\Box\rho$.
\end{prp}

\begin{prf}
For $b=0$ this is just the classical case mentioned above. For an
arbitrary $B$--field $b$ we have
$$
[\mc{D}(\Psi_+\otimes\Psi_-)]^{ev,od}_b=e^{b/2}\wedge
[\mc{D}(\Psi_+\otimes\Psi_-)]^{ev,od}=e^{b/2}\wedge
d[\Psi_+\otimes\Psi_-]^{od,ev}+e^{b/2}\wedge
d^{\star}[\Psi_+\otimes\Psi_-]^{od,ev}.
$$
The first term on the right hand side equals
$$
e^{b/2}\wedge
d[\Psi_+\otimes\Psi_-]^{od,ev}=d[\Psi_+\otimes\Psi_-]^{od,ev}_b-\frac{1}{2}e^{b/2}\wedge
db\wedge [\Psi_+\otimes\Psi_-]^{od,ev}
$$
while for the second term, we have
$$
d^{\star}[\Psi_+\otimes\Psi_-]^{od,ev}=\mp\!\star\dstar
[\Psi_+\otimes\Psi_-]^{od,ev}.
$$
Proposition~\ref{selfdual} and Lemma~\ref{hat} give
\begin{eqnarray*}
e^{b/2}\wedge d^*[\Psi_+\otimes\Psi_-]^{od,ev} & = & \mp
e^{b/2}\wedge\star
d\big(e^{b/2}\wedge\sigma[\Psi_+\otimes\Psi_-]^{ev,od}_b\big)\\
& = & \mp\frac{1}{2} e^{b/2}\wedge\star
\big(db\wedge \sigma[\Psi_+\otimes\Psi_-]^{ev,od}\big)\mp\Box_{g,b}d[\Psi_+\otimes\Psi_-]^{ev,od}_b\\
& = & \phantom{\pm}\frac{1}{2} e^{b/2}\wedge\big(db\llcorner
[\Psi_+\otimes\Psi_-]^{od,ev})+d^{\Box}[\Psi_+\otimes\Psi_-]^{od,ev}_b,
\end{eqnarray*}
hence
\begin{eqnarray*}
[\mc{D}(\Psi_+\otimes\Psi_-)]^{ev,od}_b
& = &
\phantom{+}d[\Psi_+\otimes\Psi_-]^{od,ev}_b+d^{\Box}[\Psi_+\otimes\Psi_-]_b^{od,ev}\\
& &
+\frac{1}{2}e^{b/2}\wedge\big(db\llcorner[\Psi_+\otimes\Psi_-]^{od,ev}-db\wedge [\Psi_+\otimes\Psi_-]^{od,ev}\big).
\end{eqnarray*}
Similarly, we obtain
\begin{eqnarray*}
\,[\widehat{\mc{D}}(\Psi_+\otimes\Psi_-)]^{ev,od}_b 
& = &
\pm e^{b/2}\wedge\big(d[\Psi_+\otimes\Psi_-]^{od,ev}\pm\star \dstar [\Psi_+\otimes\Psi_-]^{od,ev}\big)\\
&= &
\pm d[\Psi_+\otimes\Psi_-]_b^{od,ev}\mp
d^{\Box}[\Psi_+\otimes\Psi_-]^{od,ev}_b\\
& & \mp\frac{1}{2}e^{b/2}\wedge\big(db\llcorner
[\Psi_+\otimes\Psi_-]^{od,ev}+ db\wedge [\Psi_+\otimes\Psi_-]^{od,ev}\big).
\end{eqnarray*}
\end{prf}

We now come to our main result.

\begin{thm}\label{integrability} A generalised $G_2$--structure $(M,\rho)$ is weakly
integrable with respect to $H$ and Killing number $\lambda$ if and
only if $e^{-\phi}[\Psi_+\otimes\Psi_-]_b=\rho+\Box_{g,b}\rho$
satisfies (with $T=db/2+H$)
$$
\nabla_X\Psi_{\pm}\pm\frac{1}{4}(X\llcorner T)\cdot\Psi_+=0
$$
and
$$
\big(d\phi\pm\frac{1}{2}T\pm\lambda\big)\cdot\Psi_{\pm}=0
$$
in case of an even structure and
$$
\big(d\phi\pm\frac{1}{2}T+\lambda\big)\cdot\Psi_{\pm}=0
$$
in case of an odd structure.

The structure $e^{-\phi}[\Psi_+\otimes\Psi_-]_b=\rho+\Box_{g,b}\rho$
is strongly integrable if and only if these equations hold for
$\lambda=0$.
\end{thm}

We refer to the equation involving the covariant derivative of the
spinor as the {\em generalised Killing equation} and to the equation
involving the differential of $\phi$ as the {\em dilatino equation}.
The generalised Killing equation basically states that we have two
metric connections $\nabla^{\pm}$ preserving the underlying
$G_{2\pm}$--structures whose torsion (as it is to be defined in the
next section) is {\em skew--symmetric}. The dilatino equation then
serves to identify the components of the torsion with respect to the
decomposition into irreducible $G_{2\pm}$--modules with the
additional data $d\phi$ and $\lambda$. The generalised Killing and
the dilatino equation occur in physics as solutions to the
supersymmetry variations in type II superstring theory with bosonic
background fields~\cite{gmpw04}.

\begin{prf}
Assume that $\rho$ is even and satisfies
$d\rho=-H\wedge\rho+\lambda\hat{\rho}$. The odd case is dealt with
in a similar fashion. Applying the $\Box$--operator and using
$\Box\rho=\hat{\rho}$, we obtain
$$
d^{\Box} e^{-\phi}[\Psi_+\otimes\Psi_-]^{od}_b=e^{b/2}\wedge
\big(H\llcorner [e^{-\phi}\Psi_+\otimes\Psi_-]^{od}+\lambda
[e^{-\phi}\Psi_+\otimes\Psi_-]^{ev}\big).
$$
Consequently, Proposition~\ref{dirac} implies
\begin{equation}\label{ldirac}
[\mc{D}(e^{-\phi}\Psi_+\otimes\Psi_-)]^{ev} =  T\llcorner
[e^{-\phi}\Psi_+\otimes\Psi_-]^{od}- T\wedge
[e^{-\phi}\Psi_+\otimes\Psi_-]^{od}+\lambda
[e^{-\phi}\Psi_+\otimes\Psi_-]^{ev}.
\end{equation}
As a corollary to Proposition~\ref{lmap}, we see that
$$
T\llcorner [\Psi_+\otimes\Psi_-]^{ev,od} =
\frac{1}{8}\big(-T\cdot\Psi_+\otimes\Psi_-\pm\Psi_+\otimes T\cdot\Psi_-\mp
\sum\limits_i(e_i\llcorner T)\cdot\Psi_+\otimes
e_i\cdot\Psi_-+\sum\limits_ie_i\cdot\Psi_+\otimes[e_i\llcorner
T)\cdot\Psi_-]^{od,ev}
$$
and
$$
T\wedge [\Psi_+\otimes\Psi_-]^{ev,od} =
\frac{1}{8}[T\cdot\Psi_+\otimes\Psi_-\pm\Psi_+\otimes T\cdot\Psi_-\mp
\sum\limits_i(e_i\llcorner T)\cdot\Psi_+\otimes
e_i\cdot\Psi_--\sum\limits_ie_i\cdot\Psi_+\otimes(e_i\llcorner
T)\cdot\Psi_-]^{od,ev}.
$$
Hence~(\ref{ldirac}) entails
\begin{equation}\label{inteq1}
[\mc{D}(e^{-\phi}\Psi_+\otimes\Psi_-)]^{ev}=\frac{1}{4}e^{-\phi}[-
T\cdot\Psi_+\otimes\Psi_-+\sum\limits_ie_i\cdot\Psi_+\otimes
(e_i\llcorner  T)\cdot\Psi_-+\lambda\Psi_+\otimes\Psi_-]^{ev}.
\end{equation}
Let $D$ denote the Dirac operator associated with the Clifford
bundle $(\Delta,q)$. Contraction with $q(\cdot,e_m\cdot\Psi_+)$
yields
$$
q(De^{-\phi}\Psi_+,e_m\cdot\Psi_+)\Psi_-+e^{-\phi}\nabla_{e_m}\Psi_-=
-\frac{1}{4}e^{-\phi}q(T\cdot\Psi_+,e_m\cdot\Psi_+)\Psi_-+\frac{1}{4}e^{-\phi}(e_m\llcorner
 T)\cdot\Psi_-,
$$
and therefore
$$
e^{-\phi}\nabla_{e_m}\Psi_-=-\frac{1}{4}q(4De^{-\phi}\Psi_++e^{-\phi}
T\cdot\Psi_+,e_m\cdot\Psi_+)\Psi_-+\frac{1}{4}e^{-\phi}(e_m\llcorner
 T)\cdot\Psi_-.
$$
From this expression we deduce
$$
e_m.q(\Psi_-,\Psi_-)
=-\frac{1}{2}q(4e^{\phi}De^{-\phi}\Psi_++
T\cdot\Psi_+,e_m\cdot\Psi_+)=0,
$$
since $q\big((e_m\llcorner  T)\cdot\Psi_-,\Psi_-\big)=0$.
It follows
$$
\nabla_{e_m}\Psi_--\frac{1}{4}(e_m\llcorner  T)\cdot\Psi_-=0.
$$
We derive the corresponding expression for the spinor $\Psi_+$ by
using $\widehat{\mc{D}}$ which accounts for the minus sign.

Next we turn to the dilatino equation. From
$\nabla_X\Psi_+=-\frac{1}{4}(X\llcorner T)\cdot\Psi_+$ we deduce
$$
D\Psi_+=-\frac{3}{4}T\cdot\Psi_+.
$$
On the other hand, contracting~(\ref{inteq1}) with
$q(\,\cdot\,,\Psi_-)$ yields
$$
De^{-\phi}\Psi_+=e^{-\phi}\lambda\Psi_+-\frac{1}{4}e^{-\phi}T\cdot\Psi_+.
$$
Putting the last two equations together results in the dilatino
equation for $\Psi_+$. We can perform the same calculation with
$\widehat{\mc{D}}$ instead of $\mc{D}$ to derive the dilatino
equation for $\Psi_-$.

Conversely assume that the generalised Killing and dilatino equations of an even structure hold -- the odd case,
again, is analogous. We note that
$$
d_He^{-\phi}[\Psi_+\otimes\Psi_-]^{ev}_b=e^{-\phi}\lambda[\Psi_+\otimes\Psi_-]^{od}_b
$$
is equivalent to
\begin{equation}\label{closed}
d[\Psi_+\otimes\Psi_-]^{ev}-d\phi\wedge
[\Psi_+\otimes\Psi_-]^{ev}=\lambda
[\Psi_+\otimes\Psi_-]^{od}-T\wedge [\Psi_+\otimes\Psi_-]^{ev}.
\end{equation}
Now
\begin{eqnarray*}
d[\Psi_+\otimes\Psi_-]^{ev}-d\phi\wedge [\Psi_+\otimes\Psi_-]^{ev} &
= &
\phantom{+}\frac{1}{2}[-d\phi\cdot\Psi_+\otimes\Psi_-+\Psi_+\otimes
d\phi\cdot\Psi_-]^{od}+\\
& & +\sum\limits_i e_i\wedge\nabla_{e_i}[\Psi_+\otimes\Psi_-]^{ev}\\
& = & \phantom{+}\lambda [\Psi_+\otimes\Psi_-]^{od}+\frac{1}{4}[
T\cdot\Psi_+\otimes\Psi_-+\Psi_+\otimes
T\cdot\Psi_-]^{od}+\\
& &
+\frac{1}{2}\sum\limits_i[e_i\cdot\nabla_{e_i}\Psi_+\otimes\Psi_--
\nabla_{e_i}\Psi_+\otimes e_i\cdot\Psi_-+\\
& & +e_i\cdot\Psi_+\otimes\nabla_{e_i}\Psi_--\Psi_+\otimes e_i\cdot\nabla_{e_i}\Psi_-]^{od}\\
& = &
\phantom{+}\lambda [\Psi_+\otimes\Psi_-]^{od}-\frac{1}{8}[ T\cdot\Psi_+\otimes\Psi_-+\Psi_+\otimes  T\cdot\Psi_-\\
& & -\sum\limits_i (e_i\llcorner  T)\cdot\Psi_+\otimes
e_i\cdot\Psi_--\sum\limits_i e_i\cdot\Psi_+\otimes(e_i\llcorner
T)\cdot\Psi_-]^{od}\\
& = & \phantom{-}\lambda [\Psi_+\otimes\Psi_-]^{od} -T\wedge
[\Psi_+\otimes\Psi_-]^{ev}
\end{eqnarray*}
which proves~(\ref{closed}) and thus the theorem.
\end{prf}

\begin{rmk}
The theorem holds more generally for any 3--form $T$, closed or not.
Similarly, we can introduce a 1--form $\alpha$ and consider the
twisted differential operator $d_{\alpha}$. This substitutes $dF$ by
$dF+\alpha$ in the dilatino equation.
\end{rmk}

\subsection{The Ricci curvature and a no--go theorem}\label{ricci}

In the light of the spinorial formulation of integrability, we shall
from now on always consider the twisted differential $d_T$ applied
to a $B$--field free form $e^{-\phi}[\Psi_+\otimes\Psi_-]^{ev,od}$.
We refer to the 3--form $T$ as the {\em torsion} of the generalised
structure. Generally speaking, the {\em torsion tensor}~\cite{friv02} is defined by
$$
g(\widetilde{\nabla}_XY,Z)=g(\nabla_XY,Z)+\frac{1}{2}{\rm Tor}(X,Y,Z)
$$
and measures the difference between an arbitrary metric connection
$\widetilde{\nabla}$ and the Levi--Civita connection $\nabla$. In our situation the spinors $\Psi_+$ and $\Psi_-$ induce two $G_2$--subbundles which carry metric connections $\nabla^+$ and
$\nabla^-$ such that
$$
\nabla_X^{\pm}Y=\nabla_XY\pm\frac{1}{2}T(X,Y,\cdot).
$$
It therefore makes sense to consider the broader class of geometries
defined by two linear metric connections $\nabla^{\pm}$ with skew--symmetric and closed torsion $\pm T$. This class encapsulates all the
structures we obtain by applying the variational principle and
conveniently avoids distinguishing between ``internal" torsion $db$
coming from the ubiquitous $B$--field (corresponding to the untwisted
variational problem) and ``external" torsion $H$ which might or
might not be present. Consequently, we regard $T$ to be part of the
intrinsic data of an integrable structure. To see how the torsion
encodes the geometry, we state the following proposition (see for instance~\cite{friv03}).

\begin{prp}\label{g2tor}
For any $G_2$--structure with stable form $\varphi$ there exist
unique differential forms $\lambda\in\Omega^0(M)$,
$\theta\in\Omega^1(M)$, $\xi\in\Omega^2_{14}(M,\varphi)$, and
$\tau\in\Omega^3_{27}(M,\varphi)$ so that the differentials of
$\varphi$ and $\star\varphi$ are given by
$$
\begin{array}{ccl}
d\varphi & = & -\lambda\star\varphi+\frac{3}{4}\theta\wedge\varphi+\star\tau\\
\dstar\varphi & = &
\phantom{-}\theta\wedge\star\varphi+\xi\wedge\varphi.
\end{array}
$$
Here, $\Omega^2_{14}(M,\varphi)$ and $\Omega^3_{27}(M,\varphi)$
denote the the bundles associated with the irreducible $G_2$--modules
$\Lambda^2_{14}\leqslant\Lambda^2$ and
$\Lambda^3_{27}\leqslant\Lambda^3$~{\rm\cite{br87}}.
\end{prp}

To specify the torsion tensor of a connection is in general not
sufficient to guarantee its uniqueness. However, this is true for
$G_2$--connections with skew--symmetric torsion. Using the notation of
the previous proposition, we can assert the following result.

\begin{prp}\label{t_comp}$\!\!\!${\rm~\cite{friv02},~\cite{friv03}}\hspace{2pt}
For a $G_2$--structure with stable form $\varphi$ the following
statements are equivalent:

{\rm (i)} the $G_2$--structure is {\em integrable}, i.e. $\xi=0$.

{\rm (ii)} there exists a unique linear connection
$\widetilde{\nabla}$ whose torsion tensor ${\rm Tor}$ is skew and which
preserves the $G_2$--structure, i.e. $\widetilde{\nabla}\varphi=0$.

The torsion can be expressed by
\begin{equation}\label{tor}
{\rm Tor}=-\frac{1}{6}\lambda\cdot\varphi+\frac{1}{4}\star(\theta\wedge\varphi)-\star\tau.
\end{equation}
Moreover, the Clifford action of the torsion 3--form on the induced
spinor $\Psi$ is
\begin{equation}\label{comp}
{\rm Tor}\cdot\Psi=\frac{7}{6}\lambda\Psi-\theta\cdot\Psi.
\end{equation}
\end{prp}

Using additional subscripts $\pm$ to indicate the torsion forms of
$\nabla^{\pm}$, equations~(\ref{tor}) and~(\ref{comp}) read
$$
{\rm Tor}_{\pm}=\pm
T=-\frac{1}{6}\lambda_{\pm}\cdot\varphi_{\pm}+\frac{1}{4}\star(\theta_{\pm}\wedge\varphi_{\pm})-\star\tau_{\pm}
$$
and
$$
{\rm Tor}_{\pm}\cdot\Psi_{\pm}=\pm
T\cdot\Psi_{\pm}=\frac{7}{6}\lambda_{\pm}\Psi_{\pm}-\theta_{\pm}\cdot\Psi_{\pm}.
$$
In view of Theorem~\ref{integrability} we can use the dilatino
equation to relate the torsion components to the additional
parameters $d\phi$ and $\lambda$. We have
$$
\pm T\cdot\Psi_{\pm}=\mp2\lambda\Psi_{\pm}-2d\phi\cdot\Psi_{\pm}
$$
if the structure is even and
$$
\pm T\cdot\Psi_{\pm}=-2\lambda\Psi_{\pm}-2d\phi\cdot\Psi_{\pm}
$$
if the structure is odd.

\begin{cor}\label{torsioncomp}
If the generalised structure is weakly integrable, then there exist
two linear connections $\nabla^{\pm}$ preserving the
$G_{2\pm}$--structure with skew torsion $\pm T$. These connections
are uniquely determined. If the structure is weakly integrable and
of even type, then
$$
\begin{array}{ccl}
d\varphi_+ & = & \frac{12}{7}\lambda\star\varphi_++\frac{3}{2}d\phi\wedge\varphi_+-\star T_{27+}\\
\dstar\varphi_+ & = & 2d\phi\wedge\star\varphi_+
\end{array}
$$
and
$$
\begin{array}{ccl}
d\varphi_- & = & -\frac{12}{7}\lambda\star\varphi_-+\frac{3}{2}d\phi\wedge\varphi_-+\star T_{27-}\\
\dstar\varphi_- & = & 2d\phi\wedge\star\varphi_-,
\end{array}
$$
where $T_{27\pm}$ denotes the projection of $T$ onto
$\Omega^3_{27}(M,\varphi_{\pm})$. Moreover, the torsion can be
expressed by the formula
\begin{equation}\label{tor1}
{\rm Tor}_{\pm}=\pm T=-e^{2\phi}\star
de^{-2\phi}\varphi_{\pm}\pm2\lambda\cdot\varphi_{\pm}.
\end{equation}
If the structure is weakly integrable and of odd type, then
$$
\begin{array}{ccl}
d\varphi_+ & = & \frac{12}{7}\lambda\star\varphi_++\frac{3}{2}d\phi\wedge\varphi_+-\star T_{27+}\\
\dstar\varphi_+ & = & 2d\phi\wedge\star\varphi_+
\end{array}
$$
and
$$
\begin{array}{ccl}
d\varphi_- & = & \frac{12}{7}\lambda\star\varphi_-+\frac{3}{2}d\phi\wedge\varphi_-+\star T_{27-}\\
\dstar\varphi_- & = & 2d\phi\wedge\star\varphi_-.
\end{array}
$$
The torsion can be expressed by the formula
\begin{equation}\label{tor2}
{\rm Tor}_{\pm}=\pm T=-e^{2\phi}\star de^{-2\phi}\varphi_{\pm}
+2\lambda\cdot\varphi_{\pm}.
\end{equation}

we obtain the formulae for strongly integrable case if we set
$\lambda=0$.

Conversely, if we are given two $G_2$--structures defined by the
stable forms $\varphi_+$ and $\varphi_-$ inducing the same metric, a
constant $\lambda$ and a function $\phi$ such that~{\rm (\ref{tor1})} or~{\rm (\ref{tor2})} defines a closed 3--form $T$, then the corresponding
spinors $\Psi_{\pm}$ satisfy the integrability condition of Theorem
\ref{integrability} and hence define a (weakly) integrable
generalised $G_2$--structure (of even or odd type).
\end{cor}

The previous discussion has a striking consequence. Assume $M$ to be compact and endowed with a weakly integrable structure of even type. Then~(\ref{tor1}) and Stokes' Theorem imply
$$
\int_Me^{-2\phi}T\wedge \star T=\mp\int_MT\wedge
de^{-2\phi}\varphi_{\pm}+2\lambda\int_MT\wedge
e^{-2\phi}\star\varphi_{\pm}=\frac{4}{7}\lambda^2\int_Me^{-2\phi}vol_M.
$$
Here we have used that $dT=0$ and that the projection of
$T$ on $\varphi_{\pm}$ is given by
$T_{1{\pm}}=2\lambda\varphi_{\pm}$. The same identity holds for odd
structures. Since the left hand side is strictly positive unless
$T\equiv 0$, we need $\lambda\not=0$. As a result, we obtain the
following no--go theorem, generalising a similar statement
in~\cite{gmpw04}.

\begin{thm}\label{vanishingg2}
If $M$ is compact and carries an integrable generalised
$G_2$--structure, then $T=0$ if and only if $\lambda=0$. In this case
the spinors $\Psi_{\pm}$ are parallel with respect to the
Levi--Civita connection. Moreover, a weakly integrable structure has
necessarily non--trivial torsion.
\end{thm}

\begin{rmk}
In case of strong integrability, it follows in particular that the
underlying topological generalised $G_2$--structure cannot be exotic.
If the spinors were to be linearly dependent at one point, covariant
constancy would imply global linear dependency and we would have an
ordinary manifold of holonomy $G_2$. If the two spinors are linearly
independent at some (and hence at any) point, then the holonomy
reduces to an $SU(3)$--principal fibre bundle which is the
intersection of the two $G_2$--structures. In this case, $M$ is
locally isometric to $CY^3\times S^1$ where $CY^3$ is a Calabi--Yau
3--fold.
\end{rmk}

Next we compute the Ricci tensor. To begin with, let $\ric$ and
$\ric^{\pm}$ denote the Ricci tensors associated with the Levi--Civita
connection and the connections $\nabla^{\pm}$. Generally speaking,
the Ricci tensor $\widetilde{\ric}$ associated with a $G$--invariant
spinor $\Psi$ and a $G$--preserving, metric linear connection
$\widetilde{\nabla}$ with skew torsion is determined by the
following relation.

\begin{prp}$\!\!\!${\rm~\cite{friv02}}\hspace{2pt}
The Ricci tensor associated with $\widetilde{\nabla}$ is determined
by
$$
\widetilde{\ric}(X)\cdot\Psi=\widetilde{\nabla}_X{\rm Tor}\cdot\Psi+\frac{1}{2}(X\llcorner
d{\rm Tor})\cdot\psi
$$
and relates to the metric Ricci tensor through
$$
\ric(X,Y)=\widetilde{\ric}(X,Y)+\frac{1}{2}d^*{\rm Tor}(X,Y)+\frac{1}{4}g(X\llcorner
{\rm Tor},Y\llcorner {\rm Tor}).
$$
\end{prp}

\begin{thm}\label{riccitensor}
The Ricci--tensor of an integrable generalised $G_2$--structure is
given by
$$
\ric(X,Y)=-2H^{\phi}(X,Y)+\frac{1}{4}g(X\llcorner T,
Y\llcorner T),
$$
where $H^{\phi}(X,Y)=X.Y.\phi-\nabla_XY.\phi$ is the Hessian of the
dilaton $\phi$. It follows that the scalar curvature of the metric
$g$ is
$$
{\rm Scal}={\rm Tr}(\ric)=2\Delta
\phi+\frac{3}{4}\norm{T}^2,
$$
where $\Delta(\cdot)=-{\rm Tr_g}H^{(\cdot)}$ is the Riemannian
Laplacian.
\end{thm}

\begin{prf}
According to the previous proposition we obtain
$$
\ric(X,Y)=\frac{1}{2}\big(\ric^+(X,Y)+\ric^-(X,Y)\big)+\frac{1}{4}g(X\llcorner
T,Y\llcorner T)
$$
and it remains to compute the Ricci--tensors $\ric^+$ and $\ric^-$. Since $\nabla^{\pm}$ preserves the $G_{2\pm}$--structure, we derive 
$$
\pm(\nabla_X^{\pm}T)\cdot\Psi_{\pm}=-2(\nabla_X^{\pm}d\phi)\cdot\Psi_{\pm}
$$
from the dilatino equation, and therefore
$$
\ric^{\pm}(X)=-2\nabla_X^{\pm}d\phi.
$$
Now consider a frame that satisfies $\nabla_{e_i}e_j=0$ at a fixed
point, or equivalently,
$\nabla_{e_i}^{\pm}e_j=\pm\frac{1}{2}\sum_kT_{ijk}e_k$. As the connections $\nabla^{\pm}$ are metric, we obtain
\begin{eqnarray*}
\ric^{\pm}(e_i,e_j) & = &
-2g(\nabla_{e_i}^{\pm}d\phi,e_j)\\
& = & -2e_i.g(d\phi,e_j)+g(d\phi,\nabla^{\pm}_{e_i}e_j)\\
& = & -2e_i.e_j.\phi\pm\sum_kT_{ijk}e_k.\phi.
\end{eqnarray*}
The Hessian evaluated in this basis is just
$H^{\phi}(e_i,e_j)=e_i.e_j.\phi$, whence the result.
\end{prf}

If the dilaton is constant, then $T_{7\pm=0}$ and the underlying
$G_{2\pm}$--structures are co--calibrated, i.e.
$\dstar\varphi_{\pm}=0$. Consequently, the Ricci tensors
$\ric^{\pm}$ of $\nabla^{\pm}$ vanish. We then appeal to Theorem 5.4
of~\cite{friv02} which translated to our context asserts that if
$T_{1\pm}$ vanishes, the Levi--Civita connection reduces to the
underlying $G_{2\pm}$--structure.

\begin{prp}
The metric $g$ of a strongly integrable generalised $G_2$--structure
is Ricci--flat if and only if the dilaton $\phi$ is constant. In
particular, any strongly integrable generalised $G_2$--structure
which is homogeneous is Ricci--flat.
\end{prp}

The following example illustrates how restrictive the assumption
$dT=0$ really is. It defines a compact generalised $G_2$--structure
with constant dilaton and non--trivial $T$ such that $d_T\rho=0$,
$d_T\hat{\rho}=0$, but $dT\not=0$.

\begin{ex}
Consider the 6--dimensional nilmanifold $N$ associated with the Lie
algebra $\mf{g}$ spanned by the orthonormal basis $e_2,\ldots,e_7$
and determined by the relations
$$
de_i=\left\{\begin{array}{cl}\phantom{-}e_{37},&\,i=4\\-e_{35},&\,i=6\\\phantom{-}0,&\,\mbox{else.}\end{array}\right.
$$
We let $M=N\times S^1$ and endow $M$ with the product metric
$g=g_N+dt\otimes dt$. On $N$ we choose the $SU(3)$--structure coming
from
$$
\omega=-e_{23}-e_{45}+e_{67},\quad\psi_+=e_{356}-e_{347}-e_{257}-e_{246}
$$
which induces a generalised $G_2$--structure on $M$ with
$\alpha=e_1=dt$ (cf. Section~\ref{geng2stru}). We put
$T=e_{167}+e_{145}$. Obviously, $dT\not=0$ holds. Writing
$\rho=\omega+\psi_+\wedge\alpha-\omega^3/6$ and
$\hat{\rho}=\alpha-\psi_--\omega^2\wedge\alpha/2$ the equations
$d_T\rho=0$ and $d_T\hat{\rho}=0$ are equivalent to
$$
d\omega=0,\quad d\psi_+\wedge\alpha=-T\wedge\omega,\quad
T\wedge\alpha=d\psi_-,\quad T\wedge\psi_-=0.
$$
By design, $T\wedge\alpha=0$ and $T\wedge\psi_-=0$. Moreover,
$$
d\omega=-de_4\wedge e_5+de_6\wedge e_7=0
$$
and
$$
d\psi_-=-e_3\wedge de_4\wedge e_6+e_{34}\wedge de_6+e_{25}\wedge
de_6+e_2\wedge de_4\wedge e_7=0.
$$
Finally, we have
$$
d\psi_+\wedge\alpha=-e_{12367}-e_{12345}=-T\wedge\omega.
$$
\end{ex}

\section{Generalised $Spin(7)$--structures}\label{spin7theory}

In the same vein as generalised $G_2$--structures we can develop a theory of generalised $Spin(7)$--structures associated with
$Spin(7)\times Spin(7)$ in $\R^*\times Spin(8,8)$. We content ourselves with an
outline of this theory. Mutatis mutandis the proofs translate
without too much difficulty from the generalised $G_2$--case -- the only new feature to take into account is the chirality of the spinors.
A more detailed exposition can be found in \cite{wi04t}.

As before, a generalised $Spin(7)$--structure gives rise to two unit
spinors $\Psi_+$ and $\Psi_-$ in the associated irreducible
$Spin(8)$--representations $\Delta_+$ and $\Delta_-$ (the subscript does not indicate the chirality of the spinors). Tensoring those induces a
$Spin(7)\times Spin(7)$--invariant spinor which is given by either an
even or an odd form according to the chirality of the spinors
$\Psi_+$ and $\Psi_-$. This leads to the notion of generalised
$Spin(7)$--structures of {\em even} or {\em odd} type for spinors of equal or opposite chirality. The
box--operator is now a map
$\Box_{g,b}:\Lambda^{ev,od}\to\Lambda^{ev,od}$ for which the
(anti--)self--duality property
$$
\Box_{g,b}[\Psi_+\otimes\Psi_-]_b=(-1)^{ev,od}\Box_{g,b}[\Psi_+\otimes\Psi_-]_b
$$
holds. To any $Spin(7)\times Spin(7)$--invariant spinor we can then
associate the volume form $q(\Box_{\rho}\rho,\rho)$ and define the
dilaton $\phi$ by
$\mc{Q}(\rho)=q(\rho,\Box_{\rho}\rho)=16e^{-2\phi}vol_g$. Note that
such a spinor is not stable.

\begin{thm}
A {\em topological generalised $Spin(7)$--structure} over an
8--manifold $M$ is a topological $Spin(7)\times Spin(7)$--reduction of
the $\R^*\times Spin(8,8)$--principal fibre bundle associated with $T\oplus T^*$.
It is said to be of {\em even} or {\em odd type} according to the
parity of the $Spin(7)\times Spin(7)$--invariant spinor $\rho$. We
will denote this structure by the pair $(M,\rho)$ and call $\rho$
the {\em structure form}. Equivalently, such a structure can be
characterised by the following data:
\begin{itemize}
\item an orientation
\item a metric $g$
\item a 2--form $b$
\item a scalar function $\phi$
\item two unit half spinors $\Psi_+,\Psi_-\in\Delta$ of either equal (even type) or opposite chirality (odd type) such that
$$
e^{-\phi}[\Psi_+\otimes\Psi_-]_b=\rho.
$$
\end{itemize}
\end{thm}

The normal form description is computed as in the $G_2$--case.

\begin{prp}
There exists an orthonormal basis $e_1,\ldots,e_8$ such that
\begin{eqnarray*}
[\Psi_+\otimes\Psi_-] & = &
\phantom{+}\cos(a)+\sin(a)(e_{12}-e_{34}-e_{56}+e_{78})+\\
& &
+\cos(a)(e_{1234}+e_{1246}-e_{1278}+e_{1357}+e_{1368}+e_{1458}-e_{1467}-\\
& &
-e_{2358}+e_{2378}+e_{2457}+e_{2468}-e_{3467}+e_{4578}+e_{5678})+\\
& &
+\sin(a)(e_{1358}-e_{1367}-e_{1457}-e_{1468}+e_{2357}+e_{2368}+e_{2458}-e_{2467})+\\
& &
+\sin(a)(-e_{123456}+e_{123478}+e_{125678}-e_{345678})+\cos(a)e_{12345678},
\end{eqnarray*}
where $a=\sphericalangle(\psi_+,\psi_-)$ and the spinors are of equal chirality (even case). In the odd case, we have
\begin{eqnarray*}
[\Psi_+\otimes\Psi_-] & = & -e_1+
e_{234}+e_{256}-e_{278}+e_{357}+e_{368}+e_{458}-e_{467}+\\
& &
{}+e_{12358}-e_{12367}-e_{12457}-e_{12468}+e_{13456}-e_{13478}-e_{15678}+e_{2345678}.
\end{eqnarray*}
If the spinors are not parallel, we can express the homogeneous
components of $[\Psi_+\otimes\Psi_-]$ in terms of the invariant forms of the
intersection $Spin(7)_+\cap
Spin(7)_-=SU(4)$ (even case) or $Spin(7)_+\cap
Spin(7)_-=G_2$ (odd case). We find
$$
[\Psi_+\otimes\Psi_-] =
c+s\omega+c(\psi_+-\frac{1}{2}\omega^2)-s\psi_-\frac{s}{6}\omega^3+cvol_g
$$
where $\omega$ is the K\"ahler form and $\psi_{\pm}$
the real and imaginary parts of the holomorphic volume form stabilised by $SU(4)$. In the odd case, we have
$$
[\Psi_+\otimes\Psi_-] =
-\alpha+\varphi-\alpha\wedge\psi+\frac{1}{7}\varphi\wedge\psi
$$
where $\alpha$ denotes the dual of the unit vector in $T$ which is
stabilised by $G_2$, $\varphi$ the invariant stable 3-form on the complement and $\psi=\star_7\varphi$.
\end{prp}

Existence of generalised $Spin(7)$--structures follows from the existence of either a unit spinor (even structures) or a unit spinor and a unit vector field (odd structure), both of which is classical~\cite{lami89}.

\begin{prp}\hfill\newline
{\rm (i)} An 8--manifold $M$ carries an even topological generalised
$Spin(7)$--structure if and only if $M$ is spin and
$8\chi(M)+p_1(M)^2-4p_2(M)=0$.

{\rm (ii)} A differentiable 8--manifold $M$ carries an odd
topological generalised $Spin(7)$--structure if and only if $M$ is
spin, has vanishing Euler class and satisfies $p_1(M)^2=4p_2(M)$.
\end{prp}

For a generalised $Spin(7)$--structure of even type we can discuss
classification issues as in Section~\ref{topgeng2}. The sphere
bundle associated with $\Delta_+$ has fibre isomorphic to $S^7$ and
an 8--dimensional base, so that two transverse sections will
intersect in a curve. We meet the first obstruction for the
existence of a vertical homotopy in $H^7(M,\Z)$, which by Poincar\'e
duality trivially vanishes if $H_1(M,\Z)=0$ (e.g. if $M$ is simply
connected). Since $\pi_8(S^7)=\Z_2$, the second obstruction lies in the top cohomology module
$$
H^8(M,\pi_8(S^7))\cong\left\{\begin{array}{ll} \Z_2, & \mbox{if $M$ compact}\\
0, & \mbox{if $M$ non--compact}\end{array}\right..
$$

\begin{rmk}
The stable homotopy group $\pi_{n+k}(S^n)$ is isomorphic to the
framed cobordism group of $k$--manifolds. It is conceivable that the
$\Z_2$--class is the framed (or spin) cobordism class of the
1--manifold where the two sections coincide.
\end{rmk}

\begin{ex}
The tangent bundle of the 8--sphere is stably trivial and therefore
all the Pontrjagin classes vanish. Since the Euler class is
non--trivial, there exists no generalised $Spin(7)$--structure on
$S^8$. However, they do exist on manifolds of the form $M=S^1\times
N^7$ for $N^7$ spinnable. For instance, take $N^7$ to be the 7--sphere $S^7$. Then the tangent bundle of $M$ is trivial and
so is the sphere bundle $\mathbb{S}$ associated with the spinor
bundle. Hence $Gen(M)=[S^1\times S^7,S^7]$ which contains the set
$[S^7,S^7]=\pi_7(S^7)=\Z$. Choosing a non--trivial homotopy class in
$\pi_7(S^7)$ which we extend trivially to $S^1\times S^7$ defines an
exotic generalised $Spin(7)$--structure.
\end{ex}

Since the structure form is not stable, we cannot setup a
variational problem. We therefore follow the analogy of the
classical case and impose ad--hoc the strong integrability condition.

\begin{definition}
A topological generalised $Spin(7)$--structure $(M,\rho)$ is said to
be {\em integrable with respect to a closed 3--form} $H$ if and only
if
$$
d_H\rho=0.
$$
\end{definition}

\begin{ex}\hfill\newline
(i) Consider an 8--manifold $M$ endowed with a $Spin(7)$--invariant
4--form $\Phi$ and associated Riemannian volume form $vol_g$. The
structure form of the induced straight structure (necessarily of
even type) is given by the $B$--field transform of
$$
\rho=e^{-\phi}e^{b/2}\wedge(1-\Phi+vol),
$$
which follows from the normal form description above. We want to solve
$$
d_Te^{\phi}(1-\Omega+vol_g)=0
$$
which is equivalent to $d\phi=0,\:T=0\mbox{ and }d\Omega=0$.
Consequently, the holonomy of $(M,g)$ reduces to $Spin(7)$.

(ii) Examples of odd type are provided by product manifolds of the
form $M=N^7\times S^1$, where $N^7$ carries a $G_2$--structure
induced by the stable 3--form $\varphi$. Let $T=\eta+\xi\wedge dt$ be a
closed 4--form on $M$. The spinor
$$
\rho=e^{-\phi}e^b\wedge\big(dt\wedge(-1+\star_N\varphi)-\varphi+vol_N\big)
$$
defines a generalised $Spin(7)$--structure of odd type which is
$d_T$--closed if and only if $d\phi=0$, $T=0$ and $d\varphi=0$,
$d\star_N\varphi=-\xi\wedge\varphi$, i.e. the $G_2$--structure on $N$
is calibrated.

(iii) Manifolds with holonomy contained in $Spin(7)$ can be easily
built out of a trivial $S^1$--bundle over a 7--manifold with holonomy
$G_2$. This easily generalises to our context where a strongly
integrable generalised $G_2$--structure $(M,\rho,T)$ induces an
integrable generalised $Spin(7)$--structure of even type $(M^7\times
S^1,dt\wedge\hat{\rho}+\rho,T)$.
\end{ex}

\begin{thm}\label{spin7integrability}
A generalised $Spin(7)$--structure $(M,\rho)$ is integrable if and
only if $e^{-\phi}[\Psi_+\otimes\Psi_-]_b=\rho$ satisfies (with
$T=db/2+H$)
$$
\nabla_X\Psi_{\pm}\pm\frac{1}{4}(X\llcorner
T)\cdot\Psi_{\pm}=0
$$
and
$$
(d\phi\pm \frac{1}{2}T)\cdot\Psi_{\pm}=0.
$$
\end{thm}

Using Theorem~\ref{spin7integrability} we can discuss the torsion of
the underlying $Spin(7)$--structures defined by the invariant 4--forms
$\Omega_{\pm}$ in the same way as in the $G_2$--case. Using results
from~\cite{iv01}, we obtain:

\begin{prp}
If the generalised $Spin(7)$--structure is integrable, then
\begin{equation}\label{torspin7cor}
\pm T=e^{2\phi}\star d(e^{-2\phi}\Omega_{\pm})
\end{equation}
and
\begin{equation}\label{domegaspin7}
d\Omega_{\pm}=\frac{12}{7}d\phi\pm\star T_{48\pm},
\end{equation}
where $T_{48\pm}$ denotes the projection of $T$ onto
$\Omega^3_{48}(M,\Omega_{\pm})$, the bundle associated with the
irreducible $Spin(7)_{\pm}$--representation space $\Lambda^3_{48}$ in
$\Lambda^3$~{\rm\cite{br87}}.

Conversely, if we are given two $Spin(7)_{\pm}$--invariant forms
inducing the same metric $g$, a function $\phi$ and a closed 3--form
$T$ such that~{\rm (\ref{torspin7cor})} and~{\rm (\ref{domegaspin7})} hold, then
the corresponding spinors $\Psi_{\pm}$ satisfy
Theorem~\ref{spin7integrability} and hence define an integrable
generalised $Spin(7)$--structure.
\end{prp}

As in the $G_2$--case, we then deduce the following results.

\begin{cor}\label{vanishingspin7}
If $M$ is compact and carries an integrable generalised
$Spin(7)$--structure, then $T=0$. Consequently, the spinors
$\Psi_{\pm}$ are parallel with respect to the Levi--Civita
connection.
\end{cor}

\begin{prp}
The Ricci--tensor and the scalar curvature of an integrable
generalised $Spin(7)$--structure are given by
$$
\ric(X,Y)=-2H^{\phi}(X,Y)+\frac{1}{4}g(X\llcorner T, Y\llcorner T)
$$
and
$$
{\rm Scal}={\rm Tr}(\ric)=2\Delta \phi+\frac{3}{4}\norm{T}^2.
$$
\end{prp}

\begin{prp}
An integrable generalised $Spin(7)$--structure is
Ricci--flat if and only if the dilaton $\phi$ is constant. In
particular, any homogeneous integrable generalised
$Spin(7)$--structure is Ricci--flat.
\end{prp}

Again the condition $dT=0$ is crucial.

\begin{ex}
In conjunction with $(iii)$ of the previous example, the compact
$G_2$--structure discussed at the end of Section~\ref{ricci} gives
trivially rise to an instance of a compact generalised
$Spin(7)$--structure of even type that satisfies $d_T\rho=0$, but
$dT\not=0$. For an odd example, just take a compact calibrated
$G_2$--manifold which, for instance, can be built out of the
nilmanifold $N$ considered above. It follows that
$\dstar\varphi=\xi\wedge\varphi$ (Proposition~\ref{g2tor}), so that
$d_T\rho=0$ for $T=-\xi\wedge dt$ and $\rho$ defined as in $(ii)$ of
the previous example. 
\end{ex}

\begin{cor}
The torsion 3--form of a compact calibrated $G_2$--manifold can never be closed.
\end{cor}

\section{$T$--duality}\label{t-duality}

Type IIA and IIB string theory are interrelated by the so--called
$T$--duality. Formally speaking, $T$--duality transforms the data
$(g,b,\phi)$ consisting of a generalised metric $(g,b)$ and the
dilaton $\phi$, all living on a principal $S^1$--bundle $P\to M$ with
connection form $\theta$, into a generalised metric $(g^t,b^t)$ and
a dilaton $\phi^t$ over a {\em new} principal $S^1$--fibre bundle
$P^t\to M$ with connection form $\theta^t$. In the physics literature, the coordinate description of this dualising procedure is known under the name of
{\em Buscher rules}. A neat mathematical formulation particularly apt for applications in the setting of generalised geometry was given in~\cite{bem03}.

Let $X$ denote the dual vertical vector field of $\theta$, i.e.
$X\llcorner\theta=1$. Consider its curvature 2--form $\mc{F}$ which
we regard as a 2--form on $M$ so $d\theta=p^*\mc{F}$. We assume
to be given a closed, integral and $S^1$--invariant 3--form $T$ such
that the 2--form $\mc{F}^t$ defined by $p^*\mc{F}^t=-X\llcorner T$ is
also integral. In practice, we will assume $T=0$ so that this
condition is automatically fulfilled. Integrality of $\mc{F}^t$
ensures the existence of another principal $S^1$--bundle $P^t$, the
{\em $T$--dual} of $P$ defined by the choice of a connection form
$\theta^t$ with $d\theta^t=p^*\mc{F}^t$. Writing $T=\theta\wedge
\mc{F}^t-\mc{T}$ for a 3--form $\mc{T}\in\Omega^3(M)$, we define the
$T$--dual of $T$ by
$$
T^t=-\theta^t\wedge\mc{F}+\mc{T}.
$$
Here and from now on, we ease notation and drop the pull--back $p^*$.

To make contact with our situation, consider an $S^1$--invariant
structure form $\rho$ which we decompose into
$$
\rho=\theta\wedge\rho_0+\rho_1.
$$
The {\em $T$--dual of} $\rho$ is defined to be
$$
\rho^t=\theta^t\wedge\rho_1+\rho_0.
$$
In particular, $T$--duality reverses the parity of forms and maps
even to odd and odd to even forms. It is enacted by multiplication
with the element $X\oplus\theta\in Pin(n,n)$ on $\rho$ followed by
the substitution $\theta\to\theta^t$.

The crucial feature of $T$--duality is that it preserves the
$Spin(n,n)$--orbit structure on $\Lambda^{ev,od}TP^*$. To see this,
we decompose
$$
TP\oplus T^*\!P\cong TM\oplus\R X\oplus T^*M\oplus\R\theta,\quad
TP\oplus T^*\!P\cong TM\oplus\R X^t\oplus T^*M\oplus\R\theta^t,
$$
where $\R X$ denotes the vertical summand $\ker\pi_*$ of $TP$ which is spanned by $X$, and similarly for $\theta$ and their $T$--duals. Then consider the map $\tau:TP\oplus T^*\!P\to TP^t\oplus T^*\!P^t$
defined with respect to this splitting by
$$
\tau(V+uX\oplus\xi +v\theta)=-V+vX^t\oplus-\xi+u\theta^t.
$$
It satisfies
$$
(a\bullet\rho)^t=\tau(a)\bullet\rho^t
$$
for any $a\in TP\oplus T^*\!P$ and in particular,
$\tau(a)^2=-(a,a)$. Hence this map extends to an isomorphism
$\cliff(TP\oplus T^*\!P)\cong\cliff(TP^t\oplus T^*\!P^t)$, and any
orbit of the form $Spin(TP\oplus T^*\!P)/G$ gets mapped to an
equivalent orbit $Spin(TP^t\oplus T^*\!P^t)/G^t$ where $G$ and $G^t$
are isomorphic as abstract groups.

As an illustration of this, consider a generalised $G_2$--structure
over $P$ with structure form $\rho=\theta\wedge\rho_0+\rho_1$ and
companion $\hat{\rho}=\theta\wedge\hat{\rho}_0+\hat{\rho}_1$. These
have $T$--duals $\rho^t$ and $\hat{\rho}^t$. Since $\rho$ and
$\hat{\rho}$ have the same stabiliser $G$ inside $Spin(TP\oplus
T^*\!P)\cong Spin(7,7)$, it follows that $\rho^t$ and $\hat{\rho}^t$
are stabilised by the same $G^t$ inside $Spin(TP^t\oplus
T^*\!P^t)\cong Spin(7,7)$, which is isomorphic to $G_2\times G_2$. By
invariance, $\hat{\rho}^t$ and $\widehat{\rho^t}$ coincide up to a
constant which we henceforth ignore. The integrability condition
transforms as follows:

\begin{prp}
$$
d_{T}\rho=\lambda\hat{\rho}\mbox{ if and only if
}d_{T^t}\rho^t=-\lambda\widehat{\rho^t}
$$
\end{prp}

The proof is a straightforward computation using the definition of
$\rho$ and $\rho^t$.

We put this machinery into action as follows. Start with a
non--trivial principal $S^1$--fibre bundle $(P,\theta)$ which admits a
metric of holonomy $G_2$ or $Spin(7)$ and let $T=0$. The resulting
straight structure is strongly integrable and so is its dual, but
according to the $T$--duality rules, we acquire non--trivial torsion
given by $T^t=-\theta^t\wedge\mc{F}$. Local examples of such
$G_2$--structures exist in abundance~\cite{apsa04}. In conjunction with $(iii)$ of the first example in
Section~\ref{spin7theory}, this gives an $S^1$--invariant generalised $Spin(7)$--structure of even type with integral, $S^1$--invariant torsion $T$. Contracting with $dt$ yields $\mc{F}^t=0$, hence the $T$--dual defines an integrable generalised $Spin(7)$--structure of odd type with $T^t=T$.

As a further application of this formalism, we note that a principal
$S^1$--fibre bundle $\pi:P\to M$ whose holonomy reduces to an
$S^1$--invariant $G_2$-- or $Spin(7)$--structure must be trivial if the
base $M$ is compact. Indeed, if $\theta$ is a connection form on $P$
and $T=0$, then the $T$--dual defines a strongly integrable
generalised structure with torsion $T^t=-\theta^t\wedge\mc{F}$. But
this vanishes as a consequence of Corollaries~\ref{vanishingg2}
and~\ref{vanishingspin7}, hence $\mc{F}=0$ which implies the
triviality of $P$.

\begin{cor}
If a compact, simply--connected 7-- or a 8--manifold admits an $S^1$--invariant $G_2$-- or $Spin(7)$--structure to which the holonomy reduces, then the principal $S^1$--fibre bundle is trivial.
\end{cor}

\end{document}